\documentclass[a4paper,leqno]{article}
\usepackage{latexsym,amssymb,amsmath,amsthm}
\title{\bf Singular Homology of Arithmetic Schemes}
\author{Alexander Schmidt}
\date{July 29, 2007}

\theoremstyle{plain}
\newtheorem{theorem}{Theorem}[section]
\newtheorem{prop}[theorem]{Proposition}
\newtheorem{corol}[theorem]{Corollary}
\newtheorem{lemma}[theorem]{Lemma}

\theoremstyle{definition}
\newtheorem{defi}[theorem]{Definition}

\theoremstyle{remark}

\newtheorem{examples}[theorem]{Examples}
\newtheorem{remarks}[theorem]{Remarks}
\newtheorem{remark}[theorem]{Remark}

\def\noproof{{\unskip\nobreak\hfill\penalty50\hskip2em\hbox{}%
     \nobreak\hfill$\square$\parfillskip=0pt%
     \finalhyphendemerits=0\par}}
\def\enddemo{\ifmmode\eqno\square\else\noproof\vskip0.8truecm\fi}

\def\diagram{\renewcommand\arraystretch{1.5} $$ \begin{array}}
\def\enddiagram{\end{array} $$ \renewcommand\arraystretch{1}}

\def\rmnewname#1{\expandafter\gdef\csname#1\endcsname{{\mathop{\rm
#1}\nolimits}}}
\def\itnewname#1{{\expandafter\gdef\csname#1\endcsname{{\mathop{\it
#1}\nolimits}}}}

 \itnewname{ab}
 \itnewname{can}
 \itnewname{coker}
 \rmnewname{codim}
 \rmnewname{cone}
 \itnewname{cor}
 \rmnewname{CH}
 \rmnewname{div}
 \rmnewname{Div}
 \itnewname{et}
 \rmnewname{Ext}
 \rmnewname{Frob}
 \rmnewname{Hom}
 \itnewname{id}
 \rmnewname{Ker}
 \rmnewname{length}
 \rmnewname{Mor}
 \rmnewname{Pic}
 \rmnewname{ERes}
 \rmnewname{Quot}
 \itnewname{rec}
 \itnewname{res}
 \rmnewname{Sch}
 \rmnewname{Sm}
 \rmnewname{SmCor}
 \itnewname{sp}
 \rmnewname{Spec}
 \rmnewname{supp}
 \itnewname{tame}
 \rmnewname{top}
 \rmnewname{Tor}
 \itnewname{Zar}
 \def\A{{\mathbb A}}
 \def\C{{\underline{C}}}
 \def\F{{\mathbb F}}
 \def\G{{\mathbb G}}
 \def\ha{\underline{h}}
 
 \def\HC{{\hbox{\rm I \hspace{-6.8pt}H}}}
 \def\P{{\mathbb P}}
 \def\gP{{\mathfrak P}}
 \def\p{{\mathfrak p}}

 \def\Q{{\mathbb Q}}
 \def\R{{\hbox{\bf R}}}
 \def\Real{{\mathbb R}}
 \def\Z{{\mathbb Z}}

 \def\noi{\noindent}
 \def\bs{\bigskip}
 \def\ms{\medskip}
 \def\ds{\displaystyle}

 \def\Cal#1{{\cal #1}}

 \def\lang{\longrightarrow}
 
 \def\mapr#1{\mathrel{\mbox{$\stackrel{#1}{\longrightarrow}$}}}
 
 \def\mapd#1{\Big\downarrow\rlap{$\vcenter{\hbox{$\scriptstyle{{#1}}$}}$}}
 
 \def\liso{\mathrel{\hbox{$\longrightarrow$} \kern-12pt\lower-4pt\hbox{$\scriptstyle\sim$}\kern7pt}}
 \def\r@iso{\mathrel{\lower2pt\hbox{$\scriptstyle\sim$} \kern-8pt\hbox{$\rightarrow$}}}
 \def\riso#1{\mathrel{\stackrel{\!\!#1}{\r@iso}}}
 \def\lr@iso{\mathrel{\kern6pt\lower2pt\hbox{$\scriptstyle\sim$} \kern-12pt\hbox{$\longrightarrow$}}}
 \def\lriso#1{\mathrel{\mathop{\lr@iso}\limits^{#1}}}
 \def\eqd{\Big\|}
 
 \def\mapu#1{\Big\uparrow\rlap {$\vcenter{\hbox{$\scriptstyle{{#1}}$}}$}}
 \def\mapd#1{\Big\downarrow\rlap {$\vcenter{\hbox{$\scriptstyle{{#1}}$}}$}}
 \def\surjr#1{\stackrel{#1}{\hbox{$\relbar \! \! \! \twoheadrightarrow$}}}
 
 \def\injr#1{\stackrel{#1}{\lhook\joinrel\relbar\joinrel\rightarrow}}
 
 \def\surjd#1{\lower4pt\hbox{$\downarrow$}\kern-5.7pt\Big\downarrow\rlap {$\vcenter{\hbox{$\scriptstyle{{#1}}$}}$}}
 \def\surju#1{\lower-4pt\hbox{$\uparrow$}\kern-6.95pt\Big\uparrow\rlap {$\vcenter{\hbox{$\scriptstyle{{#1}}$}}$}}

 \def\hang{\hangindent\Itemindent}
 \def\textindent#1{\hskip\Itemindent\llap{\hbox{\rm #1}\enspace}\ignorespaces}
 \def\Item{\par\noindent\hang\textindent}
 
 \newdimen\Itemindent \Itemindent=.9cm

\begin{document}
\maketitle

\begin{minipage}{11cm}
{\footnotesize{\sc\noindent Abstract}: We construct a singular
homology theory on the category of schemes of finite type over a
Dedekind domain and verify several basic properties. For arithmetic schemes we construct a reciprocity isomorphism between the integral singular homology in degree zero and the abelianized modified tame fundamental group.}
\end{minipage}

\vskip1cm

\section{Introduction}

\indent

The objective of this paper is to construct a reasonable singular homology
theory on the category of schemes of finite type over a Dedekind domain.
Our main criterion for \lq reasonable\rq\ was to find a theory which
satisfies the usual properties of a singular homology theory and which has
the additional property that, for  schemes of finite type over
$\Spec(\Z)$, the group $h_0$ serves as the source of  a reciprocity map for tame class field theory. In the case of schemes of finite type over finite fields this role was taken over by Suslin's singular homology, see \cite{S-S1}. In this article we motivate and give
the definition of the singular homology groups  of schemes of finite type
over a Dedekind domain and we verify basic properties, e.g.\
homotopy invariance. Then we present an application to tame class field theory.

\bs
The (integral) singular homology groups $h_*(X)$ of a scheme of finite type over a
field $k$ were defined by A. Suslin as the homology of the complex
$C_*(X)$ whose $n$-th term is given by
\[
C_n(X) = \hbox{group of finite correspondences }\Delta^n_k \lang X,
\]
where $\Delta^n_k= \Spec (k[t_0,\ldots,t_n]/\sum t_i=1)$ is the
$n$-dimensional standard simplex over $k$ and a finite correspondence is a
finite linear combination $\sum n_i Z_i$ where each $Z_i$ is an integral
subscheme of $X \times \Delta^n_k$ such that the projection $Z_i \to
\Delta_k^n$ is finite and surjective. The differential $d: C_{n}(X) \to
C_{n-1}(X)$ is defined as the alternating sum of the homomorphisms which
are induced by the cycle theoretic intersection with the $1$-codimensional
faces $X \times \Delta^{n-1}_k$ in $X \times \Delta_k^{n}$. This
definition, see \cite{S-V1}, was motivated by the theorem of Dold-Thom in
algebraic topology. If $X$ is an integral scheme of finite type over the
field $\mathbb C$ of complex numbers, then Suslin and Voevodsky show that
there exists a natural isomorphism
\[
h_*(X,\Z/n\Z)\cong H_*^{sing}(X(\mathbb C),\Z/n\Z)
\]
between the algebraic singular homology of $X$ with finite coefficients
and the topological singular homology of the space $X(\mathbb C)$. If $X$
is proper and of dimension $d$, singular homology is related to the higher
Chow groups of Bloch \cite{B} by the formula $h_i(X)=\CH^d(X,i)$, see
\cite{V1}. A sheafified version of the above definition leads to the
``triangulated category of motivic complexes'', see \cite{V1}, which,
mainly due to the work of Voevodsky, Suslin and Friedlander has become a
powerful categorical framework for motivic (co)homology theories.

\ms If the field $k$ is finite and if $X$ is an open subscheme of a projective smooth variety over $k$, then we have the following relation to class field theory: there
exists a natural reciprocity homomorphism
\[
\rec: h_0(X) \lang \pi_1^t(X)^{ab}
\]
from the $0$th singular homology group to the abelianized tame fundamental
group of $X$. The homomorphism $\rec$ is injective and has a uniquely
divisible cokernel (see \cite{S-S1}, or Theorem~\ref{s-spiess} below for a more precise statement).

\ms
This connection to class field theory was the main motivation of the
author to study singular homology of schemes of finite type over Dedekind
domains. Let $S=\Spec(A)$ be the spectrum of a Dedekind domain and let $X$
be a scheme of finite type over $S$. The naive definition of singular
homology as the homology of the complex whose $n$-th term is the group of
finite correspondences $\Delta_S^n \to X$ is certainly not the correct
one. For example, according to this definition, we would have $h_*(U)=0$
for any open subscheme $U \subsetneqq S$. Philosophically, a ``standard
$n$-simplex'' should have dimension $n$ but $\Delta_S^n$ is a scheme of
dimension $(n+1)$.

\ms
If the Dedekind domain $A$ is finitely generated over a field, then one can define the
homology of $X$ as its homology regarded as a scheme over this field.

The striking analogy between number fields and function fields in one
variable over finite fields, as it is known from number theory, led to the
philosophy that it should be possible to consider {\it any} Dedekind
domain $A$, i.e.\ also if it is of mixed characteristic, as a curve over a
mysterious ``ground field'' $\F (A)$. In the case $A=\Z$ this ``field'' is
sometimes called the ``field with one element'' $\F_1$. A more precise
formulation of this idea making the philosophy into real mathematics and,
in particular, a reasonable intersection theory on ``$\Spec(\Z
\otimes_{\F_1} \Z)$'' would be of high arithmetic interest. With respect
to singular homology, this philosophy predicts that, for a scheme $X$ of
finite type over $\Spec(A)$, the groups $h_*(X)$ should be the homology
groups of a complex whose $n$-th term is given as the group of finite
correspondences $\Delta^n_{\F(A)} \to X$. Unfortunately, we do not have a
good definition of the category of schemes over $\F(A)$. To overcome this,
let us take a closer look on the situation of schemes of finite type over
a field.

\ms
Let $k$ be a field, $C$ a smooth proper curve over $k$ and let $X$ be any
scheme of finite type over $k$ together with a morphism $p:X \to C$.
Consider the complex $C_*(X;C)$ whose $n$-th term is given as

\[
\begin{array}{ccl}
 C_n(X;C)&=& \hbox{free abelian group over closed integral subschemes}\\
 && Z \subset X \times \Delta^n_k = X \times_C \Delta^n_C \hbox{ such that
  the restriction}\\
 && \hbox{of the projection }
 X \times_C \Delta^n_C \rightarrow \Delta^n_C \hbox{ to } Z \hbox{ induces a}\\
 && \hbox{finite morphism }  Z\to T \subset \Delta^n_C \hbox{ onto a closed inte-}\\
 &&  \hbox{gral subscheme } T \hbox{ of codimension 1 in }\Delta^n_C\hbox{ which}\\
 && \hbox{intersects all faces }\Delta_C^m \subset \Delta_C^n \hbox{
 properly.}
 \end{array}
\]
Then we have a natural inclusion
\[
C_*(X) \injr{} C_*(X;C)
\]
and the definition of $C_*(X;C)$ only involves the morphism $p:X \to C$
but not the knowledge of $k$. Moreover, if $X$ is affine, then both
complexes coincide.

\vskip.8cm So, in the general case, having no theory of schemes over
``$\F(A)$'' at hand, we use the above complex in order to define singular
homology. With the case $S=\Spec(\Z)$ as the main application in mind, we
define the singular homology of a scheme of finite type over the spectrum
$S$ of a Dedekind domain as the homology $h_*(X;S)$ of the complex
$C_*(X;S)$ whose $n$-th term is given by
\[
\begin{array}{ccl}
 C_n(X;S)&=& \hbox{free abelian group over closed integral
 subschemes}\\
&& Z \subset X \times_S \Delta^n_S \hbox{ such that the restriction of
the}\\
 &&  \hbox{projection }
 X \times_S \Delta^n_S \rightarrow \Delta^n_S \hbox{ to } Z \hbox{ induces a finite}\\
 && \hbox{morphism }  Z\to T \subset \Delta^n_S \hbox{ onto a closed integral}\\
 &&  \hbox{subscheme } T \hbox{ of codimension 1 in }\Delta^n_S \hbox{ which}\\
 && \hbox{intersects all faces }\Delta_S^m \subset \Delta_S^n \hbox{
 properly.}
 \end{array}
\]
The objective of this paper is to collect evidence that the so-defined
groups $h_*(X;S)$ establish a reasonable homology theory on the category
of schemes of finite type over $S$.

\ms The groups $h_*(X;S)$ are
covariantly functorial with respect to scheme morphisms and, on the
category of smooth schemes over $S$, they are functorial with respect to
finite correspondences.   If the structural morphism $p: X \to S$ factors
through a closed point $P$ of $S$, then our singular homology coincides
with Suslin's singular homology of $X$ considered as a scheme over the
field $k(P)$.

In section \ref{sect3}, we calculate the singular homology $h_*(X;S)$ if
$X$ is regular and of (absolute) dimension $1$. The result is similar to
that for smooth curves over fields. Let $\bar X$ be a regular
compactification of\/ $X$ over $S$ and\/ $Y=\bar X -X$. Then
\[ h_i(X;S)\cong H^{1-i}_{\Zar}(\bar X,\G_{\bar X,Y})
,
\]
where $\G_{\bar X,Y}=\ker (\G_{m,\bar X}\to i_* \G_{m,Y})$.

\ms
In section \ref{homotsect}, we investigate homotopy invariance. We show
that the natural projection $X \times_S \A^1_S \to X$ induces an
isomorphism on singular homology. We also show that the bivariant singular
homology groups $h_*(X,Y;S)$ (see section \ref{prelim} for their
definition) are homotopy invariant with respect to the second variable.

In section \ref{h0sect}, we give an alternative characterization of the
group $h_0$ which implies that, if $X$ is proper over $S$, we have a
natural isomorphism
\[
h_0(X;S) \cong \CH_0(X),
\]
where $\CH_0(X)$ is the group of zero-cycles on $X$ modulo rational
equivalence. Furthermore, we can verify the exactness of at least a small
part of the expected Mayer-Vietoris sequence associated to a Zariski-open
cover of a scheme $X$.

For a proper, smooth (regular?) scheme $X$ of absolute dimension $d$ over
the spectrum $S$ of a Dedekind domain,  singular homology should be
related to motivic cohomology, defined for example by \cite{V2}, by the
formula
\[
h_i(X;S) \cong H_{\Cal Mot}^{2d-i}(X, \Z(d)).
\]
For schemes over a field $k$, this formula has been proven by Voevodsky
under the assumption that $k$ admits resolution of singularities. In the
situation of schemes over the spectrum of a Dedekind domain it is true if
$X$ is of dimension~$1$ (cf.\ section \ref{sect3}). For a general $X$ it
should follow from the fact that each among the following complex
homomorphisms is a quasiisomorphism. The occurring complexes are in each
degree the free group over a certain set of cycles and we only write down
this set of cycles and also omit the necessary intersection conditions
with faces.
\diagram{c}
 C_*(X;S) \\
 \mapd{(1)} \\
  \begin{minipage}{11cm}\footnotesize ($Z \subset X \times \A^d\times \Delta^n$ projects
 finitely onto a codimension 1 subscheme in $\A^d\times
 \Delta^n$)
 \end{minipage} \\
 \mapu{(2)}\\
  \begin{minipage}
 {11cm}\footnotesize ($Z \subset X \times \A^d\times \Delta^n$ projects
 finitely onto a codimension 1 subscheme $T\subset \A^d\times
 \Delta^n$\\ such that the projection $T \to \Delta^n$ is equidimensional of relative
 dimension (d-1))
 \end{minipage}\\
 \mapd{(3)}\\
  \begin{minipage}{11cm}\footnotesize ($Z \subset X \times \A^d\times \Delta^n$
 equidimensional of relative dimension $(d-1)$ over  $\Delta^n$)
  \end{minipage}\\
  \mapu{(4)}\\
 \begin{minipage}{11cm}\footnotesize ($Z \subset X \times \A^d\times \Delta^n$
 projects quasifinite and dominant to $X \times \Delta^n$)
  \end{minipage}\\
  \mapd{(5)}\\
  \HC_{\Cal Mot}(X, \Z(d)[2d])
\enddiagram
It follows from the homotopy invariance of the bivariant singular homology
groups in the second variable, proven in  section \ref{homotsect}, that
$(1)$ is a quasiisomorphism. The statement that the other occurring
homomorphisms are also quasiisomorphisms is completely hypothetical at the
moment. However, it is, at least partly, suggested by the proof of the
corresponding formula over fields, see \cite{V1}, th.\ 4.3.7. and
\cite{F-V}, th.\ 7.1, 7.4.

\bs
We give the following application of singular homology to higher dimensional class field theory. Let $X$ be a regular connected scheme, flat and of finite type over $\Spec(\Z)$. Sending a closed point of $x$ of $X$ to its Frobenius automorphism $\Frob_x \in \pi_1^\et(X)^\ab$, we obtain a homomorphism
\[
r: Z_0(X) \lang \pi_1^\et(X)^\ab
\]
from the group $Z_0(X)$ of zero-cycles  on $X$ to the abelianized \'{e}tale fundamental group $\pi_1^\et(X)^\ab$. The homomorphism $r$ is known to have dense image. Assume for simplicity that the set $X({\mathbb R})$ of real-valued points of $X$ is empty. If $X$ is proper, then~$r$ factors through rational equivalence, defining a reciprocity homomorphism $\rec: \CH_0(X)\lang \pi_1^\et(X)^\ab$. The main result of the so-called unramified class field theory for arithmetic schemes of Bloch and  Kato/Saito \cite{K-S}, \cite{Sa} states that $\rec$ is an isomorphism
of finite abelian groups.

If $X$ is not proper, $r$ no longer factors through rational equivalence. However, consider the composite map
\[
r': Z_0(X) \mapr{r} \pi_1^\et(X)^\ab \surjr{} \pi_1^t(X)^\ab,
\]
where $\pi_1^t(X)^\ab$ is the quotient of $\pi_1^\et(X)^\ab$ which classifies finite \'{e}tale coverings of $X$ with at most tame ramification ``along the boundary of a compactification'' (see section \ref{tamesect}). We show that $r'$ factors through $h_0(X)=h_0(X;\Spec(\Z))$, defining an isomorphisms
\[
\rec: h_0(X)\liso \pi_1^t(X)^\ab
\]
of finite abelian groups. Hence the singular homology group $h_0(X)$ takes over the role of $\CH_0(X)$
if the scheme $X$ is not proper.

\bs This article was motivated by the work of A. Suslin, V. Voevodsky and
E.M. Friedlander on algebraic cycle theories for varieties over fields.
The principal ideas underlying this paper originate from discussions with
Michael Spie{\ss} during the preparation of our article \cite{S-S1}. The
analogy between number fields and function fields in one variable over
finite fields predicted that there should be a connection between the, yet
to be defined, singular homology groups of a scheme of finite type over
$\Spec(\Z)$ and its tame fundamental group, similar to that we had proven
for varieties over finite fields.  The author wants to thank M.~Spie{\ss} for fruitful discussions and for his remarks on a preliminary version of this paper.

\medskip
The bulk of this article was part of the author's Habilitationschrift at  Heidelberg University, seven years ago. However, I could not decide on publishing the material before the envisaged application to class field theory was established. This is the case now.

\section{Preliminaries} \label{prelim}

Throughout this article we consider the category $\Sch(S)$ of separated
schemes of finite type over a regular connected and Noetherian base scheme
$S$. Quite early, we will restrict to the case that $S$ is the spectrum of
a Dedekind domain, which is the main case of our arithmetic application.
We write $X \times Y= X \times_S Y$ for the fibre product of schemes $X,Y
\in \Sch(S)$. Unless otherwise specified, all schemes will be assumed
equidimensional.

\ms
Slightly modifying the approach of \cite{Fu} \S20.1, we define the
(absolute) dimension of an integral scheme $X \in \Sch(S)$ in the
following way. Let $d$ be the Krull dimension of $S$,  $K(X)$  the field
of functions of $X$ and $T$  the closure of the image of $X$ in $S$. Then
we put
\begin{equation}\nonumber
\dim X = \hbox{trdeg}(K(X)|K(T)) - \codim_S(T) +d.
\end{equation}

\begin{examples}\label{dimex}
1. Let $S= \Spec(\Z_p)$ and consider  $X= \Spec(\Z_p[T]/pT-1)
\cong \Spec(\Q_p)$, a divisor on $\A^1_S=\Spec(\Z_p[T])$. Then $\dim X= 1$ in our terminology, while $\dim_{Krull}
X =0$.

\ms\noindent 2.  The above notion of
dimension coincides with the usual Krull dimension~if

\noindent
\quad - $S$ is the spectrum of a field,

\noindent
\quad - $S$ is the spectrum of a Dedekind domain with  infinitely many
different

\noindent
\quad \,\, prime ideals (e.g.\ the ring of integers in a number field).
\end{examples}

Note that this change in the definition of dimension does not affect the
notion of codimension. For a proof of the following lemma we refer to
\cite{Fu}, lemma 20.1.

\begin{lemma}\label{fultonlemma}
\Item{\rm  (i)\,\,} Let $U \subset X$ be a nonempty open subscheme. Then
\[
\dim X = \dim U. \]
\Item{\rm (ii) \,} Let $Y$ be a closed integral
subscheme of the integral scheme $X$ over $S$. Then
\[
\dim X = \dim Y + \codim_X(Y).
\]
\Item{\rm (iii)} If $f: X \to X'$ is a dominant morphism of integral
schemes over $S$, then
\[
\dim X = \dim X' + \hbox{\rm trdeg}(K(X)|K(X')).
\]
In particular, $\dim X' \leq \dim X$ with equality if and only if $K(X)$
is a finite extension of $K(X')$.
\end{lemma}

\noindent
Recall that a closed immersion $i: Y \lang X$ is called a {\bf regular
imbedding} of codimension $d$ if every point $y$ of $Y$ has an affine
neighbourhood $U$ in $X$ such that the ideal in $\Cal O_U$ defining $Y \cap
U$ is generated by a regular sequence of length $d$. We say that two
closed subschemes $A$ and $B$ of a scheme $X$ intersect {\bf properly} if
$$\dim W =\dim A +\dim B -\dim X$$
(or, equivalently, $\codim_X W= \codim_X A + \codim_X B$) for every
irreducible component $W$ of $A \cap B$. In particular, an empty
intersection is proper. Suppose that the immersion $A \to X$ is a regular
imbedding. Then an inductive application of Krull's principal ideal
theorem shows that every irreducible component of the intersection $A \cap
B$ has dimension greater or equal to $\dim A + \dim B - \dim X$. In this
case improper intersection means that one of the irreducible components of
the intersection has a too large dimension. If $B$ is a cycle of
codimension $1$, then the intersection is proper if and only if $B$ does
not contain an irreducible component of $A$.

\ms
The {\bf group of cycles} $Z^r(X)$ (resp. $Z_r(X)$) of a scheme $X$ is the
free abelian group generated by closed integral subschemes of $X$ of
codimension $r$ (resp.\ of dimension $r$). For a closed immersion $i:Y \to
X$, we have obvious maps $i_*: Z_r(Y) \to Z_r(X)$ for all $r$. If $i$ is a
regular imbedding, we have a pull-back map
\[
i^*: Z^r(X)' \lang Z^r(Y),
\]
where $Z^r(X)' \subset Z^r(X)$ is the subgroup generated by closed
integral subschemes of $X$ meeting $Y$ properly. The map $i^*$ is given by
\[
i^* (V)= \sum_i n_i W_i,
\]
where the $W_i$ are the irreducible components of $i^{-1}(V)=V\cap Y$ and
the $n_i$ are the intersection multiplicities. For the definition of these
multiplicities we refer to \cite{Fu}, \S 6 (or, alternatively, one can use
Serre's Tor-formula \cite{Se}).

\bs
The {\boldmath\bf standard $n$-simplex} $\Delta_n=\Delta_S^n$ over $S$ is
the closed subscheme in $\A_S^{n+1}$ defined by the equation $t_0 + \cdots
+ t_n=1$. We call the sections $v_i:S \to \Delta_S^n$ corresponding to
$t_i=1$ and $t_j=0$ for $j \neq i$ the {\bf vertices} of $\Delta_S^n$.
Each nondecreasing map $\rho: [m]=\{0,1,\ldots,m\} \lang [n]=
\{0,1,\ldots,n\}$ induces a scheme morphism
\[
\bar \rho: \Delta^m \lang \Delta^n
\]
defined by $t_i \mapsto \sum_{\rho(j)=i}t_j$. If $\rho$ is injective, we
say that $\rho(\Delta^m_S)\subset \Delta^n_S$ is a {\bf face}. If $\rho$
is surjective, $\bar \rho$ is a {\bf degeneracy}. In this way
$\Delta_S^\bullet$ becomes a cosimplicial scheme. Further note that all
faces are regular imbeddings.

\vskip.8truecm The following definition was motivated in the introduction.

\begin{defi} \label{homodef}
For $X$ in $\Sch(S)$ and $n\geq 0$, the group $C_n(X;S)$ is the free
abelian group generated by closed integral subschemes $Z$ of $X \times
\Delta^n$ such that the restriction of the canonical projection
\[
 X \times \Delta^n \rightarrow \Delta^n
\]
to $Z$ induces a finite morphism $p: Z\to T \subset \Delta^n$ onto a
closed integral subscheme $T$ of codimension $d=\dim S$ in $\Delta^n$
which intersects all faces properly. In particular, such a $Z$ is
equidimensional of dimension $n$.
\end{defi}

\begin{remarks} \label{remark1} 1. If the structural morphism $X \to S$
factors through a finite morphism $S' \to S$ with $S'$ regular, then
$C_n(X;S)=C_n(X;S')$. In particular, if $S'=\{ P \}$ is a closed point of
$S$, i.e.\ if $X$ is a scheme of finite type over $\Spec(k(P))$, then
$C_n(X;S)=C_n(X;k(P))$ is the $n$-th term of the singular complex of $X$
defined by Suslin.

\ms\noi
2. If $S$ is of dimension $1$ (and regular and connected), then a closed
integral subscheme $T$ of codimension $d=1$ in $\Delta^n_S$ intersects all
faces properly if and only if it does not contain any face. If the image
of $X$ in $S$ omits at least one closed point of $S$, then this condition
is automatically satisfied.
\end{remarks}

Let $Z$ be a closed integral subscheme of $X \times \Delta^n$ which
projects finitely and surjectively onto a closed integral subscheme $T$ of
codimension $d$ in $\Delta^n$. Assume that $T$ has proper intersection
with all faces, i.e.\ $Z$ defines an element of $C_n(X;S)$.  Let $\Delta^m
\hookrightarrow \Delta^n$ be a face. Since the projection
\[
 Z \times_{\Delta_X^n} \Delta_X^m \lang
 T\times_{\Delta^n} \Delta^m
\]
is finite, each irreducible component of $Z \cap X \times
\Delta^m$ has dimension at most~$m$. On the other hand, a face is
a regular imbedding and therefore all irreducible components of
$Z \cap X \times \Delta^m$ have exact dimension $m$ and project
finitely and surjectively onto an irreducible component of $T
\cap \Delta^m$. Thus the cycle theoretic inverse image $i^*(Z)$
is well-defined and is in $C_m(X;S)$. Furthermore, degeneracy
maps are flat, and thus we obtain a simplicial abelian group
$C_\bullet(X;S)$. We use the same notation for the associated
chain complex which (in the usual way) is constructed as follows.

 Let
\[
d^i: \Delta^{n-1} \lang \Delta^n,\; i=0,\ldots, n,
\]
be the $1$-codimensional face operators defined by setting $t_i=0$.
Then we consider the complex (concentrated in positive homological
degrees)
\[
C_\bullet(X;S), \; d_n= \sum_{i=0}^n (-1)^i(d^i)^*: C_n(X;S) \to
C_{n-1}(X;S).
\]

\begin{defi}
We call $C_\bullet(X;S)$ the singular complex of $X$. Its homology groups
(or likewise the homotopy groups of $C_\bullet(X;S)$ considered as a
simplicial abelian group)
\[
h_i(X;S)= H_i(C_\bullet(X;S))\qquad \big( =\pi_i(C_\bullet(X;S))\big)
\]
are called the (integral) singular homology groups of $X$.
\end{defi}

From remark \ref{remark1}, 1. above, we obtain the following
\begin{lemma}\label{basechange}
Assume that the structural morphism $X \to S$ factors through a finite
morphism $S' \to S$ with $S'$ regular. Then for all $i$
\[
h_i(X;S) =h_i(X;S').
\]
\end{lemma}

\begin{examples} 1. If $k$ is a field and $S= \Spec(k)$, then the above definition of
$h_i(X)$ coincides with that of the singular homology of $X$ defined by
Suslin.

\smallskip\noi
2. $C_\bullet (X;S)$ is a subcomplex of Bloch's complex $z^r(X,\bullet)$,
where $r= \dim X$, and $C_\bullet(S;S)$ coincides with the Bloch complex
$z^d(S,\bullet)$. In particular,
\[
h_i(S;S)= \CH^d(S,i),
\]
where the group on the right is the higher Chow group defined by
Bloch. Note that Bloch \cite{B} defined his higher Chow groups only for equidimensional schemes over a field, but there is no problem with extending his  construction at hand.
\end{examples}

The push-forward of cycles makes $C_\bullet(X;S)$ and thus also $h_i(X;S)$
covariantly functorial on $\Sch(S)$. Furthermore, it is contravariant
under finite flat morphisms. Given a finite flat morphism $f: X' \to X$,
we thus have induced maps $f_*: h_\bullet(X';S) \to h_\bullet(X;S)$ and
$f^*: h_\bullet(X;S) \to h_\bullet(X';S)$, which are connected by the
formula $$f_* \circ f^* = \deg (f) \cdot \id_{h_\bullet(X;S)}.$$

\vskip.8cm In addition, we introduce bivariant homology groups. Let $Y$ be
equidimensional, of finite type and flat over $S$. If $X \times Y$ is
empty, we let $C_\bullet(X,Y;S)$ be the trivial complex. Otherwise, $X
\times Y$ it is a scheme of dimension $\dim X + \dim Y -d$ (as before,
$d=\dim S$) and we consider the group $C_n(X,Y;S)$ which is the free
abelian group generated by closed integral subschemes in $X \times Y
\times \Delta^n$ such that the restriction of the canonical projection
\[
 X \times Y \times \Delta^n \rightarrow Y \times \Delta^n
\]
to $Z$ induces a finite morphism $p: Z\to T \subset Y \times \Delta^n$
onto a closed integral subscheme $T$ of codimension $d$ in $Y \times
\Delta^n$ which intersects all faces\linebreak $Y \times \Delta^m$
properly. In particular, such a $Z$ is equidimensional of dimension $\dim
Y +n -d$. Furthermore, for a closed subscheme $Y' \subset Y$, we consider
the subgroup $C_n^{Y'}(X,Y;S)\subset C_n(X,Y;S)$ which is the free abelian
group generated by closed integral subschemes of $X \times Y \times
\Delta^n$ such that the restriction of the canonical projection
\[
 X \times Y \times \Delta^n \rightarrow Y \times \Delta^n
\]
to $Z$ induces a finite morphism $p: Z\to T \subset Y \times \Delta^n$
onto a closed integral subscheme $T$ of codimension $d$ in $Y \times
\Delta^n$ which intersects all faces $Y \times \Delta^m$ and all faces $Y'
\times \Delta^m$ properly.

In the same way as before, we obtain the complex $C_\bullet(X,Y;S)$, which
contains the subcomplex $C_\bullet^{Y'}(X,Y;S)$.
\begin{defi}
We call $C_\bullet(X,Y;S)$ the bivariant singular complex and its homology
groups
\[
h_i(X,Y;S)=H_i(C_\bullet(X,Y;S))
\]
the bivariant singular homology groups.
\end{defi}

Note that $C_\bullet(X,S;S)=C_\bullet(X;S)$ and
$h_i(X,S;S)=h_i(X;S)$. By pulling back cycles, a flat morphism
$Y' \to Y$ induces a homomorphism of complexes
\[
C_\bullet(X,Y;S) \lang C_\bullet (X,Y';S).
\]
If $Y' \hookrightarrow Y$ is a regular imbedding, we get a natural
homomorphism
\[
C_\bullet^{Y'}(X,Y;S) \lang C_\bullet (X,Y';S).
\]
Consider the complex of presheaves $\C_\bullet (X;S)$ which is given on
open subschemes $U \subset S$ by
\[
U \longmapsto C_\bullet (X,U;S).
\]
This is already a complex of Zariski-sheaves on $S$.

\begin{defi}
By $\ha_i(X;S)$ we denote the cohomology sheaves of the complex
$\C_\bullet(X;S)$. Equivalently, $\ha_i(X;S)$ is the Zariski sheaf on $S$
associated to
\[
U \longmapsto h_i(X,U;S).
\]
(The sheaves $\ha_i$ play a similar role as Bloch's higher Chow sheaves
\cite{B}.)
\end{defi}

Now assume that $X$ and $Y$ are smooth over $S$. By $c(X,Y)$ we denote the
free abelian group generated by integral closed subschemes $W \subset X
\times Y$ which are finite over $X$ and surjective over a connected
component of $X$. An element in $c(X,Y)$ is called a finite correspondence
from $X$ to $Y$. If $X_1,X_2,X_3$ is a triple of smooth schemes over $S$,
then (cf.\ \cite{V1}, \S2) there exists a natural composition
$c(X_1,X_2)\times c(X_2,X_3)\to c(X_1,X_3)$. Therefore one can define a
category $\SmCor(S)$ whose objects are smooth schemes of finite type over
$S$ and morphisms are finite correspondences. The category $\Sm(S)$ of
smooth schemes of finite type over $S$ admits a natural functor to
$\SmCor(S)$ by sending a morphism to its graph.

Let $X$ and $Y$ be smooth over $S$, let $\phi \in c(X,Y)$ be a finite
correspondence and let $\psi \in C_n(X,S)$. Consider the product $X \times
Y \times \Delta^n$ and let $p_1,p_2,p_3$ be the corresponding projections.
Then the cycles $(p_1 \times p_3)^*(\psi)$ and $(p_1 \times p_2)^*(\phi)$
are in general position. Let $\psi * \phi$ be their intersection. Since
$\phi$ is finite over $X$ and $\psi$ is finite over $\Delta^n$, we can
define the cycle $\phi\circ \psi$ as $(p_2\times p_3)_*(\phi * \psi)$. The
cycle $\phi \circ \psi$ is in $C_n(Y;S)$, and so we obtain a natural
pairing $c(X,Y) \times C_\bullet(X;S) \to C_\bullet(Y;S)$. We obtain the

\begin{prop}
For schemes $X$, $Y$ that are smooth over $S$, there exist
natural pairings for all $i$
\[
c(X,Y) \otimes h_i(X;S) \lang h_i(Y;S),
\]
making singular homology into a covariant functor on the category
$\SmCor(S)$.

\end{prop}

\section{Singular Homology of Curves} \label{sect3}

\indent

 We start this section by recalling some notions and lemmas from
\cite{S-V1}. Let $X$ be a scheme and let $Y$ be a closed subscheme of $X$.
Set $U=X-Y$ and denote by $i: Y\lang X $, $j:U \lang X $ the corresponding
closed and open embeddings.

We denote by $\Pic(X,Y)$ (the relative Picard group) the group whose
elements are isomorphism classes of pairs of the form $(L,\phi)$, where
$L$ is a line bundle on $X$ and $\phi: L|_Y \cong \Cal O_Y$ is a
trivialization of $L$ over $Y$, and the operation is given by the tensor
product. There is an evident exact sequence
\begin{equation}\label{relpicseq}
\Gamma(X,\Cal O_X^\times) \lang \Gamma(Y,\Cal O_Y^\times) \lang
\Pic(X,Y)\lang \Pic(X) \lang \Pic(Y).
\end{equation}

We also use the notation $\G_X$ (or $\G_m$) for the sheaf of invertible
functions on $X$ and we write $\G_{X,Y}$ for the sheaf on $X$ which is
defined by the exact sequence
\[
0 \lang \G_{X,Y} \lang \G_X \lang i_*(\G_Y) \lang 0.
\]
By \cite{S-V1}, lemma 2.1, there are natural isomorphisms
\[
\Pic(X,Y)= H^1_\Zar (X,\G_{X,Y})=H^1_\et (X,\G_{X,Y}).
\]

Assume that $X$ is integral and denote by $K$ the field of rational
functions on $X$. A relative Cartier divisor on $X$ is a Cartier divisor
$D$ such that $\supp(D)\cap Y=\varnothing$. If $D$ is a relative divisor and
$Z=\supp(D)$, then $\Cal O_X(D)|_{X-Z}=\Cal O_{X-Z}$. Thus $D$ defines an
element in $\Pic(X,Y)$. Denoting the group of relative Cartier divisors by
$\Div(X,Y)$, we get a natural homomorphism $\Div(X,Y) \to \Pic(X,Y)$. The
image of this homomorphism consists of pairs $(L,\phi)$ such that $\phi$
admits an extension to a trivialization of $L$ over an open neighbourhood
of $Y$. In particular, this map is surjective provided that $Y$ has an
affine open neighbourhood. Furthermore, we put
\[
\begin{array}{ccl}
G&=&\{f \in K^\times: f \in \ker (\Cal O_{X,y}^\times \lang \Cal
O_{Y,y}^\times) \hbox{ for any } y \in Y \}\\
&=&\{ f \in K^\times : f \hbox{ is defined and equal to $1$ at
each point of } Y \},
\end{array}
\]

The following lemmas are straightforward (cf.\ \cite{S-V1}, 2.3,2.4,2.5).

\begin{lemma}\label{i23}
Assume that $Y$ has an affine open neighbourhood in $X$. Then the
following sequence is exact:
\[
0 \lang \Gamma(X,\G_{X,Y}) \lang G \lang \Div(X,Y) \lang
\Pic(X,Y) \lang 0.
\]
\end{lemma}

\begin{lemma}\label{i24}
Assume that\/ $U$ is normal and every closed integral subscheme of $U$ of
codimension one which is closed in $X$ is a Cartier divisor (this happens
for example when $U$ is factorial). Then $\Div(X,Y)$ is the free abelian
group generated by closed integral  subschemes $T\subset U$ of codimension
one which are closed in $X$.
\end{lemma}

\begin{lemma}\label{p24}
Let $X$ be a scheme. Consider the  natural homomorphism
\[
p^*:\Pic(X) \lang \Pic(\A^1_X)
\]
which is induced by the projection $p: \A^1_X \to X$. If $X$ is reduced,
then $p^*$ is injective. If $X$ is normal, it is an isomorphism.
\end{lemma}

\begin{proof} Since $X$ is reduced, we have $p_* \G_{\A^1_X}=\G_X$.
Therefore the spectral sequence
\[
E_2^{ij}=H^i(X,R^j p_*\G_{\A^1_X}) \Longrightarrow H^{i+j}(\A^1_X,
\G_{\A^1_X})
\]
induces a short exact sequence
\[
0 \lang \Pic(X) \lang \Pic(\A^1_X) \to H^0(X, R^1p_*(\G_{\A^1_X})).
\]
This shows the first statement. The stalk of $R^1p_*(\G_{\A^1_X})$ at a
point $x\in X$ is the Picard group of the affine scheme $\Spec(\Cal
O_{X,x}[T])$. If $X$ is normal, then this group is trivial by \cite{B-M},
prop.5.5. This concludes the proof. \end{proof}

\begin{corol}\label{i25}
Assume that $X$ is normal and $Y$ is reduced. Then
\[
\Pic(X,Y)\cong \Pic(\A^1_X,\A^1_Y).
\]
\end{corol}

\begin{proof} Using the five-lemma, this follows from proposition
\ref{p24} and the exact sequence \ref{relpicseq}. \end{proof}

In the case that $S= \Spec(k)$ is the spectrum of a field $k$, our
singular homology coincides with that defined by Suslin. For a proof of
the next theorem see \cite{Li}.

\begin{theorem}\label{curve1}
Let $X$ be a smooth, geometrically connected curve over $k$, let $\bar X$ be a smooth
compactification of $X$ and let $Y=\bar X-X$. Then $h_i(X;k)=0$ for $i
\neq 0,1$ and
\[
\begin{array}{ccl} h_0(X;k)&=& \Pic(\bar X, Y),\\
 h_1(X;k)&=& \left\{ \begin{array}{cl}
                     0& \hbox{ if $X$ is affine,}\\
                     \,k^\times & \hbox{ if $X$ is proper}.
                     \end{array} \right.
 \end{array}
\]
\end{theorem}

\begin{corol}
Let $X$ be a smooth curve over a field $k$, $\bar X$ a smooth
compactification of\/ $X$ over $k$ and\/ $Y=\bar X -X$. Then for all $i$
\[
\renewcommand{\arraystretch}{1.5}
\begin{array}{rcl}
 h_i(X;k)&\cong & H^{1-i}_{\Zar}(\bar X, \G_{\bar X,Y})\\
 &\cong &{\mathbb H}_{\Zar}^{-i}\Big(\bar X, \cone\big(\G_{\bar X} \lang i_{Y*}(\G_Y)\big) \Big),
\end{array}
\renewcommand{\arraystretch}{1}
\]
where ${\mathbb H}_{\Zar}$ denotes Zariski hypercohomology.
\end{corol}

This corollary is a special case of a general duality theorem proven by
Voevodsky (\cite{V1}, th.4.3.7) over fields that admit resolution of
singularities.

\bs
We now consider the case that $S$ is the spectrum of a Dedekind domain,
which is the case of main interest for us. The proof of the following
theorem is parallel to the proof of theorem 3.1 of \cite{S-V1}, where the
relative singular homology of relative curves was calculated.

\begin{theorem} \label{curve2}
Assume that $S$ is the spectrum of a Dedekind domain and let $U$ be an
open subscheme of $S$. Let $Y\in \Sch(S)$ be regular and flat over $S$.
Setting $Y_U=Y \times U$, suppose that $Y-Y_U$ has an affine open
neighbourhood in $Y$. Then $h_i(U,Y;S)=0$ for $i \neq 0,1$ and
\[
\begin{array}{ccc}
 h_0(U,Y;S)&=& \Pic(Y, Y-Y_U), \\
 h_1(U,Y;S)&=& \Gamma(Y, {\G}_{Y,Y-Y_U}).
\end{array}
\]
\end{theorem}

\begin{proof} We may assume that $Y$ is connected. If $Y_U=Y$, then
$C_\bullet(U,Y;S)$ coincides with the Bloch complex $z^1(Y,\bullet)$. By
\cite{B}, Theorem 6.1 (whose proof applies without change to arbitrary regular schemes), we
have $h_i(U,Y;S)=0$ for $i \neq 0,1$ and
\[
\begin{array}{ccll}
 h_0(U,Y;S)&=&\Pic(Y)\\
 h_1(U,Y;S)&=&\Gamma(Y, \G_Y).
\end{array}
\]
Suppose that $Y_U \subsetneqq Y$. Then an integral subscheme $Z \subset
Y_U \times \Delta^n$ is in $C_n(U,Y;S)$ if and only if it is closed and of
codimension $1$ in $Y \times \Delta^n$. Since $Y$ is regular, such a $Z$
is a Cartier divisor and it automatically has proper intersection with all
faces (cf.\ Remark \ref{remark1},2.). Thus $C_n(U,Y;S)=\Div(Y,T)$ (see
Lemma \ref{i24}). Let $T=Y - Y_U$. If $V$ is an open affine
neighbourhood of $T$ in $Y$, then $V \times \Delta^n$ is an open affine
neighbourhood of $T\times \Delta^n$ in $Y \times \Delta^n$. According to
Lemma \ref{i23}, we have an exact sequence of simplicial abelian groups:
\begin{equation} \label{eqncurve2}
0 \to A_\bullet \rightarrow G_\bullet \to C_\bullet(U,Y;S) \rightarrow
\Pic(\Delta_Y^\bullet,\Delta_T^\bullet)\to 0,
\end{equation}
where
\[
G_n=\{ f \in k(\Delta_Y^n)^\times \,:\; f \hbox{ is defined and equal to
$1$ at each point of } \Delta_T^n \}
\]
and
$$
A_n= \Gamma(\Delta^n_Y,\G_{\Delta_Y^n,\Delta_T^n}).
$$
For each $n$, we have $A_n=A_0= \Gamma(Y,\G_{Y,T})$ and by Corollary
\ref{i25}, we have $\Pic(\Delta_Y^n,\Delta_T^n)= \Pic (Y,T)$. Let us
show that the simplicial abelian group $G_\bullet$ is acyclic, i.e.\
$\pi_*(G_\bullet)=0$. It suffices to check that for any  $f \in G_n$ such
that $\delta_i(f)=1$ for $i=0,\ldots,n$, there exists a $g \in G_{n+1}$
such that $\delta_i(g)=1$ for $i=0,\ldots,n$ and $\delta_{n+1}(g)=f$.
Define functions $g_i\in G_{n+1}$ for $i=1,\ldots,n$ by means of the
formula
\[
g_i=(t_{i+1}+\cdots +t_{n+1})+(t_0+\cdots +t_i)s_i(f).
\]
These functions satisfy the following equations:
\[
\delta_j(g_i)=\left\{
\begin{array}{lll}
\;1& \hbox{if}&j\neq i,i+1\\
(t_i+\cdots+t_n)+(t_0+\cdots+t_{i-1})f& \hbox{if}&j=i\\
(t_{i+1}+\cdots +t_n)+(t_0+\cdots t_i)f &\hbox{if}&j=i+1.
\end{array}
\right.
\]
In particular, $\delta_0(g_0)=1$, $\delta_{n+1}(g_n)=f$. Finally, we
set
\[
g=g_ng_{n-1}^{-1}g_{n-2}\cdots g_0^{(-1)^n}.
\]
This function satisfies the conditions we need. Evaluating the $4$-term
exact sequence (\ref{eqncurve2}) above, we obtain the statement of the
theorem.
\end{proof}

\begin{corol} \label{dual1} Assume that $S$ is the spectrum of a Dedekind domain.
Let $X$ be regular and quasifinite over $S$, $\bar X$ a regular
compactification of\/ $X$ over $S$ and\/ $Y=\bar X -X$. Then for all $i$
\[
\renewcommand{\arraystretch}{1.5}
\begin{array}{rcl}
 h_i(X;k)&\cong & H^{1-i}_{\Zar}(\bar X, \G_{\bar X,Y})\\
 &\cong &{\mathbb H}_{\Zar}^{-i}\Big(\bar X, \cone\big(\G_{\bar X}\lang i_{Y*}(\G_Y)\big) \Big),
\end{array}
\renewcommand{\arraystretch}{1}
\]
where ${\mathbb H}_{\Zar}$ denotes Zariski hypercohomology.
\end{corol}

\begin{proof} We may assume that $X$ is connected. By Zariski's main
theorem, $X$ is an open subscheme of the normalization $S'$ of $S$ in the
function field of $X$. As is well known, $S'=\bar X$ is again the spectrum
of a Dedekind domain and the projection $S' \to S$ is a finite morphism.
Therefore the result follows from Lemma \ref{basechange} and from
Theorem \ref{curve2} applied to the case $Y=S$.
\end{proof}

\begin{corol} Let $S$ be the spectrum of a Dedekind domain.
Assume that $X$ is regular and that the structural morphism $p:X \to S$ is
quasifinite. Let $\bar p: \bar X\to S$ be a regular compactification of\/
$X$ over $S$ and\/ $Y=\bar X -X$. Then there is a natural isomorphism
\[
\C_\bullet (X;S) \cong \bar p_* \: \G_{\bar X,Y} \,[\,1\,]
\]
in the derived category of complexes of Zariski-sheaves on $S$.
\end{corol}

\begin{proof} We may assume that $X$ is connected and we apply the result
of Theorem \ref{curve2} to open subschemes $Y \subset S$. Note that
$\bar X$ is the normalization of $S$ in the function field of $X$. The
stalk of $\ha_1(X;S)$ at a point $s \in S$ is the relative Picard group of
the semi-local scheme $\bar X \times_S S_s$ with respect to the finite set
of closed points not lying on $X$. A semi-local Dedekind domain is a
principal ideal domain, and the exact sequence (\ref{relpicseq}) from the
beginning of this section shows that also the corresponding relative
Picard group is trivial. Therefore, the complex of sheaves $\C_\bullet
(X;S)$ has exactly one nontrivial homology sheaf, which is placed in
homological degree~$1$ and is isomorphic to $\bar p_* \:\G_{\bar X,Y}$.
\end{proof}

Let us formulate a few results which easily follow from Theorem
\ref{curve2}. We hope that these results are (mutatis
mutandis) true for regular schemes $X$ of arbitrary dimension. We omit $S$
from the notation, writing $h_*(X)$ for $h_*(X;S)$ and $h_*(X,Y)$ for
$h_*(X,Y;S)$

\begin{theorem} \label{propleq1}
Let $S$ be the spectrum of a Dedekind domain. Assume that $X$ is regular
and quasifinite over $S$ (in particular, $\dim X =1$). Then the following
holds.

\bs
\Item{\rm (i)} $\;h_i(X)={\mathbb H}^{-i}_\Zar (S,\C_\bullet(X;S))$ for all $i$.\ms
\Item{\rm (ii)} (Local to global spectral sequence) There exists a spectral
sequence
\[
E_2^{ij}= H^{-i}_\Zar(S, \ha_j(X)) \Rightarrow h_{i+j}(X).
\]
\Item{\rm (iii)} (Mayer-Vietoris sequence) Let $X_1,X_2 \subset X$ be open
with $X= X_1 \cup X_2$. Then there is an exact sequence
\[
0 \to h_1(X_1 \cap X_2) \to h_1(X_1) \oplus h_1(X_2) \to h_1(X)
\hspace*{3cm}
\]
\[
\to h_0(X_1 \cap X_2) \to h_0(X_1) \oplus h_0(X_2) \to h_0(X) \to 0.
\]
\Item{\rm (iv)} (Mayer-Vietoris sequence with respect to the second variable)\\
 Let $U,V \subset S$ be open. Then there is an exact sequence
\[
0 \to h_1(X, U\cup V) \to h_1(X,U) \oplus h_1(X,V) \to h_1(X,U\cap V)
\hspace*{1cm}
\]
\[
\to h_0(X,U\cup V)\to h_0(X, U) \oplus h_0(X ,V) \to h_0(X,U\cap V) \to 0.
\]
\end{theorem}

\begin{proof} We may assume that $X$ is connected. Let $S'$ be the
normalization of $S$ in the function field of $X$, and we denote by $j_X :
X \to S'$ the corresponding open immersion (cf.\ the proof of Corollary
\ref{dual1}). Let, for an open subscheme $U\subset S$, $U'$ be its
pre-image in $S'$. Then
\[
h_i(X,U;S)=h_i(X,U';S'),
\]
and therefore we may assume that $S'=S$ in the proof of (iii) and (iv).
Then, by Corollary~\ref{dual1}, $h_i(X)=H^{1-i}_\Zar(S,\G_{S,S-X})$. Assertion
(iii) follows by applying the functor $\R\Gamma(S,-)$ to the exact
sequence of Zariski sheaves
\[
0 \lang \G_{S,S-X_1 \cap X_2} \lang \G_{S,S-X_1} \oplus \G_{S,S-X_2} \lang
\G_{S,S-X} \lang 0.
\]
For an open subscheme $j_U:U \lang S$, we denote the sheaf
$j_{U,!}j_U^*(\underline{\Z})$ by $\Z_U$. Then, for a sheaf $F$ on $S$, we
have a canonical isomorphism
\[
H^i_\Zar(U,j_U^*F)\cong \Ext^i_{S}(\Z_U,F).
\]
Applying the functor $\R\Hom_{S}(-,\G_{S,S-X})$
to the exact sequence of Zariski sheaves
\[
0 \lang \Z_{U\cap V} \lang \Z_U \oplus \Z_V \lang \Z_{U \cup V}
\lang 0,
\]
Theorem \ref{curve2} implies assertion (iv).
From (iv) it follows that the complex $\C_\bullet(X)$ is pseudo-flasque in
the sense of \cite{B-G}, which shows assertion (i). Finally, (ii) follows
from the corresponding hypercohomology spectral sequence converging to
${\mathbb H}^{-i}_\Zar (S,\C_\bullet(X;S))$ and from (i).
\end{proof}

 Finally, we deduce an exact Gysin sequence for one-dimensional schemes.
In order to formulate it, we need the notion of twists. Let $\G_m$ denote
the multiplicative group scheme $\A^1_S-\{0\}$ and let $X$ be any scheme
of finite type over $S$. For $i=1, \cdots, n$, let $D^i_\bullet(X \times
\G_m^{\times (n-1)};S)$ be the direct summand in $C_\bullet(X \times
\G_m^{\times n};S)$ which is given by the homomorphism
\[
\G_m^{\times (n-1)} \lang \G_m^{\times n}, \quad (x_1,\ldots,x_{n-1})
\mapsto (x_1,\ldots, 1_i, \ldots, x_{n-1})
\]
We consider the complex $C_\bullet(X \times \G_m^{\wedge n};S)$ which is
defined as the direct summand of the complex $C_\bullet(X \times
\G_m^{\times n};S)$ complementary to the direct summand $\sum_{i=1}^n
D^i_\bullet(X \times \G_m^{\times(n-1)};S)$\quad (cf.\ \cite{S-V2}, \S3).
\begin{defi} For  $n\geq0$, we put
\[
h_i(X(n);S)= H_{i+n}(C_\bullet(X \times \G_m^{\wedge n};S)).
\]
\end{defi}
In particular, we have $h_i(X(0);S)=h_i(X;S)$ for all $i$ and
$h_i(X(n);S)=0$ for $i < -n$. If $X=\{P\}$ is a closed point on $S$, then
(see \cite{S-V2}, lemma 3.2.):
\[
h_i(\{P\}(1);S)=\left\{ \begin{array}{ll}
                      k(P)^\times & \hbox{for }i=-1,\\
                      \; 0 & \hbox{otherwise.}
                      \end{array}\right.
\]
The next corollary follows from this and from Theorem \ref{curve2}.

\begin{corol}
Assume that $X$ is regular and quasifinite over $S$ and that $U$ is an
open, dense subscheme in $X$. Then we have a natural exact sequence
\[
0\to h_1(U) \to h_1(X) \to h_{-1}((X-U)(1)) \to h_0(U) \to h_0(X) \to 0.
\]
\end{corol}

\section{Homotopy Invariance} \label{homotsect}

\indent

Throughout this section we fix our base scheme $S$, which is the spectrum
of a Dedekind domain, and we omit it from the notation, writing $h_*(X)$
for $h_*(X;S)$ and $h_*(X,Y)$ for $h_*(X,Y;S)$. Our aim is to prove that
the relative singular homology groups $h_*(X,Y)$ are homotopy invariant
with respect to both variables.

\begin{theorem}\label{homot}
Let $X$ and $Y$ be of finite type over $S$. Then the projection $X \times
\A^1 \to X$ induces isomorphisms
\[
 h_i(X \times \A^1,Y) \liso h_i(X,Y)
\]
for all $i$.
\end{theorem}

Let $i_0,i_1: Y \lang Y \times \A^1$ be the embeddings defined by the
points (i.e.\ sections over $S$) $0$ and $1$ of $\A^1=\A^1_S$.

\ms
Recall that $\Delta^n$ has coordinates $(t_0,\ldots,t_n)$ with $\sum
t_i=1$. Vertices are the points (i.e.\ sections over $S$)
$p_i=(0,\ldots,0,1,0,\ldots,0)$ with $1$ in the $i$th place. Consider the
linear isomorphisms
\[
\theta_i : \Delta^{n+1} \lang \Delta^n \times \A^1, \quad i=0,\ldots,n
\]
which are defined by taking $p_j$ to $(p_j,0)$ for $j\leq i$ and to
$(p_{j-1},1)$ if $j>i$. Then consider for each $n$ the formal linear
combination
\[
T_n= \sum_{i=0}^n (-1)^i\theta_i.
\]
Let us call a subscheme $F \subset \Delta^n \times \A^1$ a face if it
corresponds to a face in $\Delta^{n+1}$ under one of the linear
isomorphisms $\theta_i$. Using this terminology, $T_n$ defines a
homomorphism from a subgroup of $C_n(X,Y \times \A^1)$ to $C_{n+1}(X,Y)$.
This subgroup is generated by cycles which have good intersection not only
with all faces $Y \times \A^1 \times \Delta^m$ but also with all faces of
the form $Y \times F$, where $F$ is a face in $\A^1\times\Delta^n$.

\bs
We will deduce Theorem \ref{homot} from the

\begin{prop} \label{i0undi1} The two chain maps
\[
i_{0*}, i_{1*}: C_\bullet(X,Y)\lang C_\bullet (X\times \A^1,Y)
\]
are homotopic. In particular, $i_{0*}, i_{1*}$ induce the same map on
homology.
\end{prop}

\begin{proof}
 Let $D \subset \A^1 \times \A^1$ be the diagonal.
Consider the map
\[
V_n: C_n(X,Y) \lang C_n(X \times \A^1, Y \times \A^1)
\]
which is defined by sending a cycle $Z \subset X \times Y \times \Delta^n$
to the cycle $Z \times D\subset X \times Y \times \Delta^n \times \A^1
\times \A^1$. If $Z$ projects finitely and surjectively onto $T \subset Y
\times \Delta^n$, then $Z \times D$ projects finitely and surjectively
onto $T \times \A^1\subset Y \times \Delta^n \times \A^1$. Therefore $V_n$
is well-defined. Fortunately, $T \times \A^1$ has proper intersection with
all faces $Y \times F$, where $F$ is a face in $\Delta^n \times \A^1$.
Therefore the composition
\[
T_{n*} \circ V_n: C_n(X,Y) \lang C_n(X\times \A^1,Y\times A^1)\lang
C_{n+1}(X \times \A^1, Y)
\]
is well-defined for every $n$. These maps give the required homotopy.
\end{proof}

\begin{proof}[Proof of Theorem \ref{homot}] Let $\tau: \A^1 \times \A^1 \lang
\A^1$ be the multiplication map. Consider the diagram
\diagram{ccc}
C_\bullet(X\times \A^1,Y)&\mapr{p_*}&C_\bullet(X,Y)\\
\mapd{i_{0*},i_{1*}}&&\mapd{i_{0*},i_{1*}}\\ C_\bullet(X \times \A^1\times
\A^1,Y)&\mapr{\tau_*}&C_\bullet(X \times \A^1,Y).
\enddiagram
We have the following equalities of maps on homology:
\[
i_{0*}\circ p_* = \tau_* \circ i_{0*} = \tau_* \circ i_{1*} =
\id_{h_\bullet(X,Y)}.
\]
Therefore, $p_*$ is injective on homology. But on the other hand, $p \circ
i_0= \id_X$, which shows that $p_*$ is surjective. This concludes the
proof.\end{proof}

Now, exploiting a moving technique of \cite{B}, we prove that the
bivariant singular homology groups $h_*(X,Y)$ are homotopy invariant with
respect to the second variable.

\begin{theorem}\label{hobase} Assume that $S$ is the spectrum of a
Dedekind domain and
let $X$ and $Y$ be of finite type over $S$. Then the projection $p: Y
\times \A^1
\to Y$ induces isomorphisms for all $i$
\[
h_i(X,Y) \liso h_i(X,Y \times \A^1).
\]
\end{theorem}

A typical intermediate step in proving a theorem like \ref{hobase} would
be to show that the induced chain maps $i_0^*, i_1^*: C_\bullet(X,Y\times
\A^1)\lang C_\bullet (X,Y)$ are homotopic. However, $i_0^*, i_1^*$ are
only defined as homomorphisms on the subcomplex:
\[
i_0^*, i_1^*: C_\bullet^{Y\times \{0,1\}}(X,Y\times \A^1)\lang C_\bullet
(X,Y).
\]
(The maps $T_n^*: C_n(X,Y\times\A^1)\lang C_{n+1}(X,Y)$ would define a
homotopy $i_0^*\sim i_1^*: C_n(X,Y\times\A^1)\lang C_{n}(X,Y)$, if all
these maps {\it would} be defined.)

\bs
The proof of Theorem \ref{hobase} will consist of several steps. First,
we show that the inclusion
\[
C_\bullet^{Y\times \{0,1\}}(X,Y\times \A^1)\lang C_\bullet(X,Y\times \A^1)
\]
is a quasiisomorphism. Then we show that the homomorphisms
\[
i_0^*, i_1^*: C_\bullet^{Y\times \{0,1\}}(X,Y\times \A^1)\lang C_\bullet
(X,Y)
\]
induce the same map on homology. Finally, we deduce Theorem \ref{hobase}
from these results.

\bigskip
In the proof we will apply a moving technique of \cite{B} which
was used there to show the homotopy invariance of the higher Chow
groups. As long as we have to deal with cycles of codimension
$1$, this technique also works in our more general situation (this
is the reason for the restriction to the case that $S$ is the
spectrum of a Dedekind domain).

\ms
We would like to construct a homotopy between the identity of the complex
$C_\bullet (X,Y\times \A^1)$ and another map which takes its image in the
subcomplex $C_\bullet^{Y\times \{0,1\}} (X,Y\times \A^1)$. What we {\it
can} do is the following:

\ms
For a suitable scheme $S'$ over $S$ we construct a homotopy between the
pullback map $C_\bullet (X,Y\times \A^1)\lang C_\bullet (X,Y\times
\A^1\times S')$ and another map whose image is contained in the subcomplex
$C_\bullet^{Y\times \{0,1\}\times S'} (X,Y\times \A^1\times S')$.
(Eventually, we will use $S'=\A^1_S$ but perhaps this would be too many
$\A^1$'s in the notation.)

\ms Let  (for the moment) $\pi: S'\to S$ be any integral scheme of finite type over $S$
and let $t$ be an element in
$\Gamma(S',\Cal O_{S'})$. Consider the action
\[
 \A^1_{S'} \times_{S'} (Y \times \A^1)_{S'} \lang (Y \times \A^1)_{S'}
\]
of the smooth group scheme $ \A^1_{S'}$ on $(Y \times \A^1)_{S'}$ given by
additive translation
\[
a \cdot (y,b)) = (y,a+b)
\]
and consider the morphism $\psi: \A^1_{S'} \to \A^1_{S'}$ given by
multiplication by $t$: $a
\mapsto ta$. The points $0,1$ of $\A^1_{S'}$ give rise to isomorphisms
\[
\psi(0),\psi(1): (Y \times \A^1)_{S'} \lang (Y \times \A^1)_{S'}
\]
($\psi(0)$ is the identity and $\psi(1)$ sends $(y,b)$ to $(y,t+b)$).
Furthermore, setting $\phi(y,a,b)= (y,\psi(b)\cdot a, b)$, we obtain an isomorphism
\[
\phi: (Y \times \A^1 \times \A^1)_{S'} \lang (Y \times \A^1 \times
\A^1)_{S'}.
\]

\bs\noi
We would like to compose the maps
\[
C_n (X,Y\times \A^1) \mapr{\pi^*} C_n(X,Y\times \A^1 \times S')
\mapr{pr^*} C_n(X,(Y\times \A^1)\times \A^1 \times S')
\]
\[ \mapr{\phi^*}
C_n(X,(Y\times \A^1)\times \A^1 \times S') \mapr{T_n^*} C_{n+1}(X,Y\times
\A^1 \times S'),
\]
but we are confronted with the problem that the map $T_n^*$ is not defined
on the whole group $C_n(X,(Y\times \A^1)\times \A^1 \times S')$. The next
proposition tells us that the composition is well-defined if $S'=\A^1_S=
\Spec\, S[t]$.

\begin{prop} \label{homotconstr} Suppose that  $S'=\A^1_S= \Spec\, S[t]$. Then the
composition
$$
H_n= T_n^* \circ \phi^* \circ pr^* \circ \pi^*: C_n (X,Y\times \A^1) \lang
C_{n+1}(X,Y\times \A^1 \times S')
$$
is well-defined for every $n$. The family $\{H_n\}_{n\geq 0}$ defines a
homotopy
\[
\pi^*=\psi(0)\circ \pi^* \sim \psi(1)\circ \pi^*: C_n (X,Y\times \A^1)
\lang C_{n}(X,Y\times \A^1 \times S').
\]
Furthermore, the image of the map $\psi(1)\circ \pi^*$ is contained in the
subcomplex\linebreak $C_{\bullet}^{Y \times \{0,1\} \times S'}(X,Y\times
\A^1 \times S')$.
\end{prop}

\begin{proof} Recall that all groups $C_\bullet$ are relative to the base
scheme $S$ which we have omitted from the notation. At the moment, the map
$H_n$ is only defined as a map to the group of cycles in $X\times Y
\times\A^1\times \Delta^{n+1} \times S'$. If $Z \subset X \times Y \times
\A^1 \times \Delta^n$ projects finitely and surjectively onto an
irreducible subscheme $T \subset Y \times \A^1 \times \Delta^n$ of
codimension one, then $\phi^*\circ pr^* \circ \pi^*(Z)$ projects finitely
and surjectively onto the irreducible subscheme of codimension one
$T'=\phi^*\circ pr^* \circ \pi^*(T) \subset (Y \times \A^1) \times
\Delta^n \times \A^1\times S'$. Therefore, in order to show that $H_n(Z)$
is in $C_{n+1}(X,Y\times \A^1 \times S')$, we have to check that
$\theta_i^{-1}(T')$ has proper intersection with all faces for
$i=0,\ldots,n$. Thus we have to show that $T'$ has proper intersection
with all faces $(Y\times \A^1)\times F\times S'$, where $F$ is a face in
$\Delta^n \times \A^1$ (as defined above). Since $T'$ has codimension one,
this comes down to show that it does not contain any irreducible component
of any face (we did not assume $Y$ to be irreducible, but we can silently
assume that it is reduced). Consider the projection
\[
Y \times \A^1 \times \Delta^n \times \A^1 \times S' \lang S'.
\]
We can check our condition
by considering the fibre over the generic point of $S'$. More
precisely,
let $k$ be the function field of $S$ and let $K=k(t)$ be the function
field of $S'$. Let $(Y_1)_k, \ldots, (Y_r)_k$ be the irreducible
components of $Y_k$. Then an irreducible subscheme $T'\subset Y \times
\A^1\times  \Delta^n \times \A^1
\times S'$ of codimension one meets all faces $Y \times \A^1 \times F
\times S'$ ($F$ a face of $\Delta^n \times \A^1$)  properly if and only if
$T_K$ does not contain $(Y_i\times \A^1)_K \times_K F_K$ for $i=1,\ldots,r$.

Now we arrived exactly at the situation considered in \cite{B},
\S2. The result follows from \cite{B}, Lemma 2.2,
by taking $(Y \times \A^1)_k$ for the scheme $X$ of
that lemma, taking $\A^1_k$ as the algebraic group $G$ acting on
$X$ by additive translation on the second factor and choosing the
map $\psi: \A^1_K\to G_K$ of that lemma as the morphism which
sends $a$ to $ta$. The fact that the $H_n$ define the homotopy is
a straightforward computation.

It remains to show that the image of the map $\psi(1)\circ \pi^*$ is
contained in the subcomplex $C_{\bullet}^{Y \times \{0,1\} \times
S'}(X,Y\times \A^1 \times S')$. But this is a again a condition which says
that a subscheme of codimension one does not contain certain subschemes.
In the same way as above, this can be verified over the generic fibre of
$S'$, and the result follows from the corresponding statement of \cite{B},
Lemma 2.2.
 \end{proof}

\begin{corol}\label{corquis1}
The natural inclusion
\[
C_\bullet^{Y\times \{0,1\}}(X,Y \times \A^1) \lang C_\bullet(X,Y \times
\A^1)
\]
is a quasiisomorphism.
\end{corol}

\begin{proof} Let $S'=\A^1_S$. Then the homomorphism

\ms\noi
$\pi^*: C_\bullet(X,Y \times \A^1)/C_\bullet^{Y\times \{0,1\}}(X,Y \times
\A^1)$
\[
 \lang
C_\bullet(X,Y \times \A^1\times S')/C_\bullet^{Y\times \{0,1\}\times
S'}(X,Y \times \A^1 \times S')
\]
is nullhomotopic (the $H_n$ of Proposition \ref{homotconstr} give the
homotopy). In order to conclude the proof, it suffices to show that the
nullhomotopic homomorphism $\pi^*$ is injective on homology. Suppose that
for a cycle $z$ in degree $n$ we have $\pi^*(z)=d_n(w)$. Then we find an
$a \in \Gamma(S,\Cal O_S)$ such that the specialization (i.e.\ $t \mapsto
a$) $w(a)$ is well-defined. But then $z=d_n(w(a))$.
\end{proof}

\begin{prop} \label{psi1}
Suppose that  $S'=\A^1_S= \Spec\, S[t]$. Then the
composition
\[
C_n^{Y \! \times \! \{0,1\}}(X,Y \! \times \! \A^1) \stackrel{\psi(1)\circ
\pi^*}{\lang} C_n^{Y \! \times \! \{0,1\} \! \times \! S'}(X,Y \! \times
\! \A^1 \! \times \! S') \mapr{T_n^*} C_{n+1}(X,Y \! \times \! S')
\]
is well-defined, giving a homotopy
\[
i_0^* \circ \psi(1)\circ \pi^* \sim i_1^* \circ \psi(1)\circ \pi^* :
C_\bullet^{Y \times \{0,1\}}(X,Y\times\A^1)\lang C_{\bullet}(X,Y \times
S').
\]
\end{prop}

\begin{proof} Let again $k$ be the function field of $S$ and let $K=k(t)$
be that of $S'$. We use the following fact, which is explained in the
proof of \cite{B}, cor.\ 2.6:

\smallskip\noi
If $z_k$ is a cycle on $Y \times \A^1\times \Delta^n$ which intersects all
faces $(Y \times \A^1\times \Delta^m)_k$ properly, then $\psi(1)\circ
\pi^*(z_k)\subset (Y \times \A^1 \times \Delta^n)_K$ intersects all faces
$(Y \times F)_K$ (where $F$ is any face in $\A^1 \times \Delta^n$)
properly.

\smallskip\noindent
We deduce the statement of Proposition \ref{psi1} from this in the same
manner as we deduced Proposition \ref{homotconstr} from \cite{B}, Lemma
2.2.  The fact that the maps $T_n^*\circ \psi(1)\circ \pi^*$ define the
homotopy is a straightforward computation.\end{proof}

\begin{corol}\label{i0i1}
The two maps
\[
i_0^*, i_1^*: C_\bullet^{Y\times \{0,1\}}(X,Y\times \A^1)\lang C_\bullet
(X,Y)
\]
induce the same map on homology.
\end{corol}

\begin{proof} Consider the commutative diagram
\diagram{ccc} C_\bullet^{Y
\times \{0,1\}} (X,Y \times \A^1)&\mapr{\pi^*}& C_\bullet^{Y \times
\{0,1\}\times S'} (X,Y \times \A^1\times S')\\
\mapd{i_0^*,i_1^*}&&\mapd{i_0^*,i_1^*}\\
C_\bullet(X,Y)&\mapr{\pi^*}&C_\bullet(X,Y \times S').
\enddiagram
The same specialization argument as in the proof of Corollary
\ref{corquis1} shows that $\pi^*$ is injective on homology. Therefore it
suffices to show that $i_0^*\circ \pi^*=i_1^*\circ \pi^*$ on homology. By
Proposition \ref{homotconstr}, we have a homotopy $\pi^* \sim
\psi(1)\circ \pi^*$, and hence it suffices to show that the maps $i_0^*
\circ \psi(1)\circ \pi^* $ and $i_1^* \circ \psi(1)\circ \pi^*$ induce the
same map on homology. But this follows from Proposition \ref{psi1}.
\end{proof}

Now we conclude the proof of Theorem \ref{hobase}. First of all, note
that
$$p^*(C_\bullet(X,Y))\subset C_\bullet^{Y \times \{0,1\}}(X,Y \times
\A^1)$$ and that $i_0^*\circ p^*=id$, such that $p^*$ is injective on
homology. Consider the multiplication map
\[
\tau: \A^1 \times \A^1 \lang \A^1.
\]
It is flat and therefore $\tau^*$ exists. Consider the diagram
\diagram{ccc} C_\bullet(X,Y \times \A^1)&\mapr{\tau^*}&C_\bullet (X,
Y\times \A^1 \times \A^1)\\ \mapu{q.iso.}&&\mapu{q.iso.}\\ C_\bullet^{Y
\times \{0,1\}}(X,Y\times \A^1)&
\stackrel{\tau^*}{\cdots\!\!\!>}&C_\bullet^{Y \times
\A^1\times\{0,1\}}(X,Y \times \A^1\times \A^1)\\
\mapd{i_0^*,i_1^*}&&\mapd{i_0^*,i_1^*}\\
C_\bullet(X,Y)&\mapr{p^*}&C_\bullet(X,Y \times \A^1).
\enddiagram

\noi
One easily observes that $\tau^*$ sends a cycle $z \in C_n^{Y \times
\{0,1\}}(X,Y \times \A^1)$ to a cycle in $C_n^{Y \times
\A^1\times\{0,1\}}(X,Y \times \A^1\times \A^1)$ and that for such a $z$
the following equalities hold:
\begin{eqnarray}\label{zyk1}
i_0^*\circ \tau^* (z) &=& p^* \circ i_0^* (z)\\ \label{zyk2} i_1^*\circ
\tau^* (z)&=&z.
\end{eqnarray}
By Corollary \ref{corquis1}, any class in $h_n(X,Y\times\A^1)$ can be
represented by an element in $C_n^{Y \times \{0,1\}}(X,Y\times \A^1)$.
Therefore (\ref{zyk1}) shows that, in order to prove that $p^*$
is surjective on homology, it suffices to show that $i_0^* \circ \tau^*$
is. But, by Corollary~\ref{i0i1}, $i_0^* \circ \tau^*$ induces the same map on
homology as $i_1^* \circ \tau^*$, which is the identity, by (\ref{zyk2}).
\enddemo

 A naive definition of homotopy between scheme morphisms
is the following: Two scheme morphisms $\phi,\psi: X \lang X'$ are
homotopic if there exists a morphism
\[
H: X \times \A^1 \lang X'
\]
with $\phi= H\circ i_0$ and $\psi= H \circ i_1$. (This is not an
equivalence relation!) The next corollary is an immediate consequence of
Proposition \ref{i0undi1}.

\begin{corol}
If two morphisms
\[
\phi,\psi: X \lang X'
\]
are homotopic, then they induce the same map on singular homology, i.e.\
for every scheme $Y$ flat and of finite type over $S$, the homomorphisms
\[
\phi_* ,\psi_*: h_i(X,Y) \lang h_i(X',Y)
\]
coincide for all $i$.
\end{corol}

\ms
Now we recall the definition of relative singular homology from
\cite{S-V1}. Suppose that $Y$ is an integral scheme and that $X$ is any
scheme over $Y$.

For $n\geq 0$, let $C_n(X/Y)$ be the free abelian group generated by
closed integral subschemes of $X \times_Y \Delta^n_Y$ such that the
restriction of the canonical projection
\[
 X \times_Y \Delta^n_Y \lang \Delta^n_Y
\]
to $Z$ induces a finite surjective morphism $p: Z\to \Delta^n_Y$. Let $i:
\Delta^m_Y \hookrightarrow \Delta^n_Y$ be a face. Then all irreducible
components of $Z \cap X \times_Y \Delta^m_Y$ have the ``right'' dimension
and thus the cycle theoretic inverse image $i^*(Z)$ is well-defined and in
$C_m(X/Y)$. Furthermore, degeneracy maps are flat, and thus we obtain a
simplicial abelian group $C_\bullet(X/Y)$.
 As above, we use the same notation for the complex of abelian groups
obtained by taking the alternating sum of face operators. The groups
\[
h_i(X/Y)= H_i(C_\bullet(X/Y))
\]
are called the {\bf relative singular homology groups} of $X$ over $Y$.

\vskip.8cm
We have seen in section \ref{prelim} that singular homology is
covariantly functorial on the category $\SmCor(X)$ of smooth schemes over
$S$ with finite correspondences as morphisms. For $X,Y \in \Sm(S)$ the
group of finite correspondences $c(X,Y)$ coincides with $C_0(X \times
Y/Y)$ and we call two finite correspondences  homotopic if they have the
same image in $h_0(X\times Y/Y)$. The next proposition shows that
homotopic finite correspondences induce the same map on singular homology.

\begin{prop} \label{homotcorr}
For smooth schemes $X,Y\in \Sm(S)$, the natural pairing $$c(X,Y) \otimes
h_i(X;S) \to h_i(Y;S)$$ factors through $h_0(X\times Y/Y)$, defining
pairings
\[
h_0(X\times Y /Y) \otimes h_i(X;S) \lang h_i(Y;S)
\]
for all \/$i$.
\end{prop}

\begin{proof} Let $W \subset X \times Y \times \Delta^1= X \times Y \times
\A^1 $ define an element in $C_1(X\times Y/Y)$. Let $W^j=i_j^*(W)$, for
$j=0,1$, so that $d_1(W)= W^0 -W^1 \in C_0(X \times Y/Y)$. Let $\psi \in
C_n(X;S)$. We have to show that $(W^0,\psi)=(W^1,\psi)$. Considering $W$
as an element in $C_0(X \times Y \times \A^1/Y \times \A^1)$, the
composite $(W,\psi)$ is in $C_n^{\{0,1\}}(Y,\A^1;S)$ and
$(W^j,\psi)=i_j^*((W,\psi))$ for $j=0,1$. Therefore the result follows
from Corollary \ref{i0i1}. \end{proof}

\section{\boldmath Alternative characterization of $h_0$} \label{h0sect}

 For a noetherian scheme $X$ we have the identification
\[
\CH^d(X,0)=\CH^d(X)
\]
between the higher Chow group $\CH^d(X,0)$ and the group $\CH^d(X)$ of
$d$-co\-dimen\-sio\-nal cycles on $X$ modulo rational equivalence (see
\cite{Na}, prop.3.1). Fixing the notation and assumptions of the previous
sections, we now give an analogous description for the group $h_0(X;S)$.

Let $C$ be an integral scheme over $S$ of absolute dimension $1$. Then to
every rational function $f\neq 0$ on $C$, we can attach the zero-cycle
$\div(f)\in C_0(C;S)$ (see \cite{Fu}, Ch.I,1.2). Let $\tilde{C}$ be the
normalization of $C$ in its field of functions. Denoting the normalization
morphism by $\phi:\tilde{C} \to C$, we have $\phi_*(\div(f))=\div(f)$. If
$C$ is regular and connected, then we denote by $P(C)$ the regular
compactification of $C$ over $S$, i.e.\ the uniquely determined regular
and connected scheme of dimension $1$ which is proper over $S$ and which
contains $C$ as an open subscheme.

With this terminology, for an integral scheme $C$ of absolute
dimension~$1$, elements in the function field $k(C)$ are in 1-1
correspondence to morphisms $P(\tilde C)\to \P^1_S$, which are not
$\equiv\infty$.

\begin{theorem} \label{h0char}
The group $h_0(X;S)$ is the quotient of the group of zero-cycles on $X$
modulo the subgroup generated by elements of the form $\div(f)$, where

\ms
{ \Item{-} $C$ is a closed integral curve on $X$, \ms
 \Item{-} $f$ is a rational function on $C$ which, considered as a
 rational function on $P(\tilde{C})$, is defined and $\equiv 1$ at every
 point of $P(\tilde{C})-\tilde{C}$.

}
\end{theorem}

\begin{proof}
 We may suppose that $X$ is reduced. Let $Z \subset X \times
\Delta^1$ be an integral curve such that the projection $Z \to \Delta^1$
induces a finite and surjective morphism of $Z$ onto a closed integral
subscheme $T$ of codimension $1$ in $\Delta^1$. Embed $\Delta^1$ linearly
to $\P^1=\P^1_S$ by sending $(0,1)$ to $0=(0:1)$ and $(1,0)$ to
$\infty=(1:0)$. Since $Z \to \Delta^1$ is finite, the projection $ Z \to
\P^1$ corresponds to a rational function $g$ on $Z$ which is defined and
$\equiv 1$ at every point of $P(\tilde Z) - Z$. Let $\bar Z$ be the
closure of $Z$ in $X \times \P^1$, and let $\bar C$ be the image of $\bar
Z$ under the (proper) projection $X \times \P^1 \to X$, considered as a
reduced (hence integral) subscheme of $X$.

\ms
We have to consider two cases:

\smallskip\noindent
1. If $\bar C=P$ is a closed point on $X$, then $Z= \{P\} \times \Delta^1$
and $d_1(Z)=0$.

\smallskip\noindent
2. If $\bar C$ is an integral curve, then the image $C$
of $Z$ under $X \times \P^1\to X$ is an open subscheme of $\bar C$.
Consider the extension of function fields
\[
k(Z)|k(C)
\]
and let $f\in k(C)$ be the norm of $g$ with respect to this extension.
Then $f$ is defined and $\equiv 1$ at every point of $P(\tilde C) -C$ and
\[
\div(f)=\delta_0(Z)-\delta_1(Z)= d(Z).
\]
If $X$ is of dimension $1$, the last equality follows from \cite{Na},
prop.1.3. The general case can be reduced to this by replacing $X$ by
$\bar C$. Considering $f$ as a rational function on $\bar C$, it satisfies
the assumption of the theorem.

\ms
It remains to show the other direction. Let $C$ and $f$ be as in the
theorem. We have to show that $\div(f)\in C_0(X;S)$ is a boundary. To see
this, interpret $f$ as a nonconstant morphism $U \to \P^1$ defined on an
open subscheme $U \subset C$ and let $\bar Z$ be the closure of the graph
of this morphism in $X \times \P^1$. The scheme $\bar Z$ is integral, of
dimension $1$ and projects birationally and properly onto $C$. Consider
again the open linear embedding $\Delta^1 \subset \P^1$ which is defined
by sending $(0,1)$ to $0$ and $(1,0)$ to $\infty$ and let $Z= \bar Z \cap
X \times \Delta^1$. The properties of $f$ imply that the induced
projection $Z \to \Delta^1$ is finite and surjective onto a closed
subscheme of codimension $1$ in $\Delta^1$, thus defining an element of
$C_1(X;S)$. Finally note that $d(Z)=\delta_0(Z)-\delta_1(Z)=\div(f)$.
\end{proof}

This immediately implies the following corollary.

\begin{corol}\label{h0ch0}
If $X$ is proper over $S$, then
\[
h_0(X;S) = \CH_0(X).
\]
\end{corol}

\begin{corol} \label{erzdD1}
The natural homomorphism
\[
\bigoplus_{i_C} d(C_1(C;S)) \mapr{{i_C}_*} d(C_1(X;S))
\]
is surjective, where $i_C: C \to X$ runs through all $S$-morphisms from a
regular scheme $C$ over $S$ of dimension $1$ to $X$.
\end{corol}

\begin{proof} By Theorem \ref{h0char}, $d(C_1(X;S)$ is generated by
elements of the form $\div(f)$, where $f$ is a rational function on an
integral curve on $X$ satisfying an additional property. The normalization
$\tilde C$ of $C$ is a regular scheme of dimension $1$ and let $i:\tilde C
\to X$ the associated morphism. Considering $f$ as a rational function on
$\tilde C$, we have the equality
\[
i_*(\div(f))= \div (f).
\]
By the additional property of $f$, the associated line bundle $\Cal L(\div
(f))$ over the compactification $P(\tilde C)$ together with its canonical
trivialization over $P(\tilde C)- \tilde C$ defines the trivial element in
$\Pic(P(\tilde C), P(\tilde C)-\tilde C)$. Therefore, the calculation of
singular homology of regular schemes of dimension $1$ (see Theorems
\ref{curve1} and \ref{curve2}), shows that $\div (f)$ is in
$d(C_1(C;S))$. This finishes the proof. \end{proof}

Now we can prove the exactness of a part of the Mayer-Vietoris sequence
for $X$ of arbitrary dimension.

\begin{prop}
Let $S=U \cup V$ be a covering by Zariski-open subschemes $U$ and $V$.
Then the natural sequence
\[
h_0(X;S)\lang h_0(X;U) \oplus h_0(X;V) \lang h_0(X; U\cap V) \lang 0
\]
is exact.
\end{prop}

\begin{proof} First of all, the homomorphism
\[
C_0(X;U)\oplus C_0(X;V) \lang C_0(X;U\cap V)
\]
is surjective, and therefore so is $h_0(X;U)\oplus  h_0(X;V) \lang
h_0(X;U\cap V)$.

\bs Consider the commutative diagram

{\footnotesize
\diagram{cccccccc}
 0&&0&&0&&\\
 \mapd{}&&\mapd{}&&\mapd{}&&\\
 d(C_1(X;S))&\lang&d(C_1(X;U))\oplus d(C_1(X;V)) &\lang&d(C_1(X;U\cap
 V))\\
 \mapd{}&&\mapd{}&&\mapd{}&&\\
 C_0(X;S)&\injr{}&C_0(X;U)\oplus C_0(X;V)&\lang&C_0(X;U\cap V)&\lang
&0\\
 \mapd{}&&\mapd{}&&\mapd{}&&\\
 h_0(X;S)&\lang&h_0(X;U)\oplus h_0(X;V)&\lang& h_0(X;U\cap V)&\lang&0\\
 \mapd{}&&\mapd{}&&\mapd{}&&\\
 0&&0&&0&&.
\enddiagram
}

\noindent
The middle row and the middle and right columns are exact. Therefore the
snake lemma shows that the lower line is exact if and only if the
homomorphism
\begin{equation} \label{dD1}
 d(C_1(X;U))\oplus d(C_1(X;V)) \lang d(C_1(X;U\cap
 V))
\end{equation}
is surjective. By Theorem \ref{propleq1},(iv), we observe that
(\ref{dD1}) is surjective if $X$ is regular and of dimension $1$. For a
general $X$, put $X'= X \times_S (U\cap V)$. Then the commutative diagram
\diagram{ccc}
 C_1(X';U)\oplus C_1(X';V) &\lang& C_1(X';U\cap V) \\
 \mapd{}&&\eqd \\
 C_1(X;U)\oplus C_1(X;V) &\lang & C_1(X;U\cap V)
\enddiagram
shows that, in order to show the surjectivity of (\ref{dD1}), we may
suppose that $X=X'$. Now the statement follows from Corollary
\ref{erzdD1}, using the commutative diagram
\diagram{ccc} \ds\bigoplus_{i_C} d(C_1(C;U)) \oplus d(C_1(C;V))&
\mapr{{i_C}_*} &d(C_1(X;U)) \oplus d(C_1(X;V))\\
 \surjd{} && \mapd{}\\
\ds\bigoplus_{i_C} d(C_1(C;U\cap V))& \surjr{{i_C}_*} & d(C_1(X;U\cap V)).
\enddiagram
This concludes the proof.\end{proof}

A similar argument shows the

\begin{prop}
Let $X=X_1 \cup X_2$ be a covering by Zariski open subschemes $X_1$ and
$X_2$. Then the natural sequence
\[
h_0(X_1\cap X_2;S)\lang h_0(X_1;S) \oplus h_0(X_2;S) \lang h_0(X; S) \lang
0
\]
is exact.
\end{prop}

\begin{proof} We omit the base scheme $S$ from our notation. First of all,
the homomorphism
\[
C_0(X_1)\oplus C_0(X_2) \lang C_0(X)
\]
is surjective, and therefore so is $h_0(X_1)\oplus h_0(X_2) \lang h_0(X)$.

\bs
Consider the commutative diagram {\footnotesize
\diagram{cccccccc}
 0&&0&&0&&\\
 \mapd{}&&\mapd{}&&\mapd{}&&\\
 d(C_1(X_1 \cap X_2))&\lang&d(C_1(X_1))\oplus d(C_1(X_2)) &\lang&d(C_1(X))\\
 \mapd{}&&\mapd{}&&\mapd{}&&\\
 C_0(X_1\cap X_2)&\injr{}&C_0(X_1)\oplus C_0(X_2)&\lang&C_0(X)&\lang
&0\\
 \mapd{}&&\mapd{}&&\mapd{}&&\\
 h_0(X_1 \cap X_2)&\lang&h_0(X_1)\oplus h_0(X_2)&\lang& h_0(X&\lang&0\\
 \mapd{}&&\mapd{}&&\mapd{}&&\\
 0&&0&&0&&.
\enddiagram
}

\noindent
The middle row and the middle and right columns are exact. Therefore the
snake lemma shows that the lower line is exact if and only if the
homomorphism
\begin{equation} \label{dD1*}
 d(C_1(X_1))\oplus d(C_1(X_2)) \lang d(C_1(X))
\end{equation}
is surjective. By Theorem \ref{propleq1},(iii), we observe that
(\ref{dD1*}) is surjective if $X$ is regular and of dimension $1$.

\bs
For a morphism $i:C\to X$ we use the notation $C_1=i^{-1}(X_1)$ and
$C_2=i^{-1}(X_2)$, thus $C=C_1 \cup C_2$ is a Zariski open covering.

Now the required statement for arbitrary $X$ follows from Corollary
\ref{erzdD1}, using the commutative diagram
\diagram{ccc} \ds\bigoplus_{i_C} d(C_1(C_1)) \oplus d(C_1(C_2))&
\mapr{{i_C}_*} &d(C_1(X_1)) \oplus d(C_1(X_2))\\
 \surjd{} && \mapd{}\\
\ds\bigoplus_{i_C} d(C_1(C))& \surjr{{i_C}_*} & d(C_1(X)).
\enddiagram
 This concludes the proof. \end{proof}

We conclude this section with the following surjectivity result.
\begin{prop} Let $X$ be regular and
let $U$ be a dense open subscheme in $X$. Then the natural homomorphism
\[
h_0(U;S) \lang h_0(X;S)
\]
is surjective.
\end{prop}

\begin{proof} Let $P$ be a $0$-dimensional point on $X$ which is
not contained in $U$. We have to show that the image of $P$ in
$h_0(X;S)$ is equal to the image of a finite linear combination
$\sum n_i P_i$ with $P_i \in U$ for all $i$. Choose a
one-dimensional subscheme $C$ on $X$ such that $P$ is a regular
point on $C$ and such that $C$ is not contained in $X-U$. We find
such a curve, since $X$ is regular: Indeed, $\Cal O_{X,P}$ is a
$d$-dimensional regular local ring, with $d=\dim X$. Let
$\mathfrak m$ be the maximal ideal and $\mathfrak a$ the ideal
defining the closed subset $(X-U)\cap \Spec(\Cal O_{X,P})$.
Choose elements $\bar x_1, \ldots, \bar x_{d-1}$ in $\mathfrak
m/\mathfrak m^2$ which span a $(d-1)$-dimensional subspace which
does not contain $\mathfrak a+ \mathfrak m/\mathfrak m$. Lifting
$\bar x_1, \ldots, \bar x_{d-1}$ to a regular sequence $x_1,
\ldots, x_{d-1}\in \Cal O_{X,P}$, the ideal $(x_1, \ldots,
x_{d-1})$ is a prime ideal of height $(d-1)$ which does not
contain $\mathfrak a$. Finally, extend $\bar x_1, \ldots,
x_{d-1}$ to an affine open neighbourhood of $P$ in $X$ and choose
$C$ as the closure of their zero-locus.

Consider the normalization $\tilde C$ of $C$ and let $P(\tilde C)$ be a
regular compactification over $S$. Let $P(\tilde C)- \tilde C =\{
P_1,\ldots,P_r\}$ and let $P_{r+1}, \ldots, P_{s}$ be the finitely many
closed points on $\tilde C$ mapping to $C\cap (X-U)$. Let $\tilde P$ be
the unique point on $\tilde C$ projecting to $ P \in C$. Let
$D=\{P_1,\ldots,P_s, \tilde P\}$ and consider the ring $A=\Cal O_{P(\tilde
C),D}$, which is a semi-local principal ideal domain. We find an element
$f \in A$ which has exact valuation $1$ at $\tilde P$ and which is $\equiv
1$ at each $P_i$, $i=1,\ldots, n$. Then $(\div f) \subset X$ is of the
form $P+ \sum Q_i$ with $Q_i \in U$.
\end{proof}

\section{Review of tame coverings}\label{tamesect}

 The concept of tame ramification stems from number theory: A
finite extension of number fields $L|K$ is called tamely ramified at a
prime $\gP$ of $L$ if the associated extension of completions
$L_\gP|K_\gP$ is a tamely ramified extension of local fields. The latter
means that the ramification index is prime to the characteristic of the
residue field. It is a classical result that composites and towers of
tamely ramified extensions are again tamely ramified. This concept
generalizes to separable extensions of arbitrary discrete valuation fields by requiring that the associated residue field extensions are separable.

\ms
Let from now on $S$ be the spectrum of an excellent Dedekind domain and let $X\in \Sch(S)$. Our aim is to say when a finite \'{e}tale covering $Y\to X$  is tame. Here ``tame'' means tamely ramified along the boundary of a compactification $\bar X$ of $X$ over $S$. If $\bar X$ is regular and $D=\bar X - X$ is a normal crossing divisor, then one can use the approach of \cite{SGA1}, \cite{G-M}:

\begin{defi}[\cite{G-M}, 2.2.2] \label{TR1}
A finite \'{e}tale covering $Y \to X$ is called tame (along $D$) if the extension of
function fields $k(Y)|k(X)$ is tamely ramified at the discrete
valuations associated to the irreducible components of $D$.
\end{defi}

Even if one restricts attention to regular schemes, one is confronted with the following problems

\begin{itemize}
\item if $X$ is regular, we do not know whether there exists a regular compactification with a NCD as boundary,
\item the notion of tameness might depend on the choice of the compactification $\bar X$ of $X$.
\item even if the first two questions can be answered in a positive way, there is no obvious functoriality for the tame fundamental group (already  for an open immersion).
\end{itemize}

All these problems are void in the case of a regular curve $C$, where a canonical compactification $\bar C$ exists. Starting from the therefore obvious notion of tame coverings of  regular curves,  G.~Wiesend \cite{W1} proposed the following definition.

\begin{defi} \label{TR3} Let $X$ be a separated integral scheme of finite type over $S$. A finite \'{e}tale covering $Y\to X$ is called {\bf tame} if for every integral curve $C \subset X$ with normalization $\tilde C\to C$ the base change
\[
Y\times_X \tilde C \longrightarrow \tilde C
\]
is a tame covering of the regular curve $\tilde C$.
\end{defi}
This definition has the advantage of making no use of a compactification of $X$. Furthermore, it is obviously stable under base change. However, it is difficult to decide whether a given \'{e}tale covering is actually tame. For coverings of normal schemes several authors (cf.\ \cite{Ab}, \cite{C-E}, \cite{S1}) made suggestions for a definition of tameness which all come down to the following notion, which we want to call numerically tameness here.

\begin{defi}
Let $\bar X\in \Sch(S)$ be  normal connected and proper, and let $X \subset \bar X$ be an open subscheme. Let $Y\to X$ be a finite \'{e}tale Galois covering and let $\bar Y$ be the normalization of $\bar X$ in the function field $k(Y)$ of $Y$. We say that $Y \to X$ is {\bf numerically tame} (along $D=\bar X -X$) if the order of the inertia group $T_x(\bar Y|\bar X) \subset Gal(\bar Y|\bar X)=Gal(Y|X)$ of each closed point $x\in D$ (see \cite{B-Comm} Ch.\,V, \S 2.2 for the definition of inertia groups) is prime to the residue characteristic of $x$. A finite \'{e}tale covering $Y \to X$ is called numerically tame if it can be dominated by a numerically tame Galois covering.
\end{defi}

\begin{prop}\label{ntimpliestame}
Let $\bar X\in \Sch(S)$ be  normal connected and proper, and let $X \subset \bar X$ be an open subscheme. If the finite \'{e}tale covering $Y\to X$ is numerically tame (along $\bar X-X$), then it is tame.
\end{prop}

\begin{proof}
For regular curves the notions of tameness and of numerically tameness obviously coincide. Therefore the statement of the proposition follows from the fact that numerically tame coverings are stable under base change, see \cite{S1}.
\end{proof}

The following theorems \ref{compare} and \ref{nilpot} are due to G.~Wiesend.

\begin{theorem}[\cite{W1}, Theorem~2] \label{compare}
Assume that $\bar X$ is regular and that $D=\bar X -X$ is a NCD. Then, for a finite \'{e}tale covering $Y\to X$, the following condition are equivalent

\ms
\Item{\rm  (i)\,\,} $Y \to X$ is tame according to Definition \ref{TR1},
\Item{\rm  (ii)\,\,} $Y \to X$ is tame (according to Definition \ref{TR3}),
\Item{\rm  (iii)\,\,} $Y \to X$ is numerically tame.
\end{theorem}

\begin{remark}
The equivalence of (i) and (iii) had already been shown in \cite{S1}.
\end{remark}

\begin{theorem}[\cite{W1}, Theorem~2] \label{nilpot}
Assume that $\bar X$ is regular (but make no assumption on $D=\bar X -X$). If a numerically tame covering $Y\to X$ can be dominated by a Galois covering with nilpotent Galois group, then it is tame.
\end{theorem}

In particular, for nilpotent coverings of a regular scheme $X$ the notion of numerically tameness does not depend on the choice of a regular compactification $\bar X$ (if it exists). This had already been shown in \cite{S1}. A counter-example with non-nilpotent Galois group can be found in \cite{W1}, Remark~3.

\section{Finiteness results for tame fundamental groups}
The tame coverings of a connected integral scheme $X\in Sch(S)$ satisfy the axioms of a Galois category (\cite{W1}, Proposition~1). After choosing a geometric point $\bar x$ of $X$ we have the fibre functor $(Y\to X)\mapsto \Mor_X(\bar x, Y)$ from the category of tame covering of $X$ to the category of sets, whose automorphisms group is called the {\bf tame fundamental group} $\pi_1^t(X,\bar x)$. It classifies finite connected tame coverings of $X$.  We have an obvious surjection
\[
\pi_1^\et(X,\bar x) \twoheadrightarrow \pi_1^t(X,\bar x),
\]
which is an isomorphism if $X$ is proper. Assume that $X$ is normal, connected and let $\bar X$ be a normal compactification.  Then, replacing tame coverings by numerically tame coverings, we obtain in an analogous way the {\bf numerically tame fundamental group} $\pi_1^{nt}(\bar X,\bar X-X, \bar x)$, which classifies finite connected numerically tame coverings of $X$ (along  $\bar X-X$).  By Proposition~\ref{ntimpliestame} we have a surjection
\[
\varphi:\ \pi_1^t(X,\bar x) \twoheadrightarrow \pi_1^{nt}(\bar X,\bar X-X, \bar x),
\]
which, by Theorem~\ref{nilpot}, induces an isomorphism on the maximal pro-nilpotent factor groups if $\bar X$ is regular. If, in addition, $\bar X - X$ is a normal crossing divisor then $\varphi$ is an isomorphism by Theorem~\ref{compare}. The fundamental groups of a connected scheme $X$ with respect to different base points are isomorphic, and the isomorphism is canonical up to inner automorphisms. Therefore, when working with the maximal abelian quotient of the \'{e}tale fundamental group (tame fundamental group, n.t.\ fundamental group) of a connected scheme, we are allowed to omit the base point from notation.

\ms
Now we specialize to the case $S=\Spec(\Z)$, i.e.\ to arithmetic schemes. In \cite{S1} we proved the finiteness of the abelianized numerically tame fundamental group  $\pi_1^{nt}(\bar X,\bar X-X)^\ab$ of a connected normal scheme, flat and of finite type over $\Spec(\Z)$ with respect to a normal compactification $\bar X$. The proof given there can be adapted to apply also to the larger group $\pi_1^t(X)^\ab$.

\begin{theorem}\label{finiteness}
Let $X$ be a connected normal scheme, flat and of finite type over $\Spec(\Z)$.  Then the abelianized  tame fundamental group $\pi_1^{t}(X)^\ab$ is finite.
\end{theorem}

For the proof we need the following two lemmas. The first one extends  \cite{S1}, Corollary 2.6 from numerical tameness to tameness.

\begin{lemma} \label{exdiv} Let $X\in Sch(S)$ be normal and connected, $p$ a prime number and  $Y\to X$ a finite \'{e}tale Galois covering whose Galois group is a finite $p$-group. Let $\bar X$ be a normal compactification of\/ $X$ and assume there exists a prime divisor $D$ on $\bar X$ which is ramified in $k(Y)|k(X)$ and which contains a closed point of residue characteristic $p$. Then $Y\to X$ is not tame.
\end{lemma}

\begin{proof}
The statement of the lemma is part of the proof of \cite{W1}, Theorem~2.
\end{proof}

\begin{lemma}\label{zpext}
Let $A$ be a strictly henselian discrete valuation ring with perfect (hence algebraically
closed) residue field and with quotient field\/ $k$. Let $k_\infty|k$ be a $\Z_p$-extension.
Let $K|k$ be a regular field extension and let $B\subset K$ be a discrete valuation ring
dominating $A$. Then $B$ is ramified in $Kk_\infty|K$.
\end{lemma}

\begin{proof}
See \cite{S1}, Lemma 3.2.
\end{proof}

\begin{proof}[Proof of Theorem \ref{finiteness}] The proof is a modification of the proof of \cite{S1}, Theorem~3.1.  Let $\bar X$ be a normal compactification of $X$ over $\Spec(\Z)$. Let $k$ be the normalization of $\Q$ in the function field of $X$ and put $S=\Spec(\Cal O_k)$. Then the natural projection $\bar X \to \Spec(\Z)$ factors through $S$.

Since $X$ is normal, for any open subscheme $V$ of
$X$ the natural homomorphism $\pi_1^\et(V)\to \pi_1^\et(X)$ is surjective. Therefore also the
homomorphism
\[
\pi_1^t(V)^{ab} \lang \pi_1^t(X)^{ab}
\]
is surjective and so we may replace $X$ by a suitable open subscheme and assume that $X$ is
smooth over $S$. Let $T\subset S$ be the image of $X$. Consider the
commutative diagram
 \diagram{ccccccc}
 0&\lang & \Ker(X/T)& \lang& \pi_1^\et(X)^{ab} & \lang &\pi_1^\et(T)^{ab}\\
 &&\mapd{}&&\surjd{} &&\surjd{}\\
 0&\lang & \Ker^t(X/T)& \lang& \pi_1^t(X)^{ab} &
 \lang &\pi_1^t(T)^{ab}
 \enddiagram
where the groups $\Ker(X/T)$ and $\Ker^t(X/T)$ are defined by the exactness of the
corresponding rows, and the two right vertical homomorphisms are surjective. By a theorem of
Katz and Lang (\cite{K-L}, th.1), the group $\Ker(X/T)$ is finite. By classical
one-dimensional class field theory, the group $\pi_1^t(T)^{ab}$ is finite (it
is the Galois group of the ray class field of $k$ with modulus $\prod_{\p \notin T}\p$). The
kernel of $\pi_1^\et(T)^{ab}\to \pi_1^t(T)^{ab}$ is generated by the ramification
groups of the primes of $S$ which are not in $T$. Denoting the product of the residue
characteristics of these primes by $N$, we see that $\pi_1^\et(T)^{ab}$ is the product of a
finite group and a topologically finitely generated pro-$N$ group. Therefore the same is also
true for $\pi_1^\et(X)^{ab}$ and for $\pi_1^t(X)^{ab}$. Hence it suffices
to show that the cokernel $C$ of the induced map $\Ker(X/T) \to \Ker^t(X/T)$ is a torsion
group.

Let $K$ be the function field of $X$ and let $k_1$ be the maximal abelian extension of $k$
such that the normalization $X_{Kk_1}$ of $X$ in the composite $Kk_1$ is ind-tame over $X$.  By \cite{K-L}, lemma 2, (2), the normalization of $T$ in $k_1$ is
ind-\'{e}tale over $T$. Let $k_2|k$ be the maximal subextension of $k_1|k$ such that the
normalization $T_{k_2}$ of $T$ in $k_2$ is tame over $T$.  Then $G(k_2|k)=\pi_1^t(T)^{ab}$ and $C\cong G(k_1|k_2)$.

In order to show that $C$ is a torsion group, we therefore have to show that $k_1|k_2$ does
not contain a $\Z_p$-extension of $k_2$ for any prime number $p$. Since $k_2|k$ is a finite
extension and $k_1|k$ is abelian, this is equivalent to the assertion that $k_1|k$ contains
no $\Z_p$-extension of $k$ for any prime number $p$.

Let $p$ be a prime number and suppose
that $k_\infty|k$ is a $\Z_p$-extension such that the normalization $X_{Kk_\infty}$ is
ind-tame over $X$. A $\Z_p$-extension of a number field is unramified outside $p$ and there exists at least one ramified prime dividing $p$, see e.g.\ \cite{NSW}, (10.3.20)(ii).  Let $k'$ be the maximal unramified subextension of $k_\infty|k$ and let $S'$ be
the normalization of $S$ in $k'$. Then the base change $\bar X'=\bar X\times_{S}S'\to X$
is \'{e}tale. Hence $\bar X'$ is normal and the pre-image $X'$ of $X$ is smooth and geometrically
connected over $k'$. So, after replacing $k$ by $k'$, we may suppose that $k_\infty|k$ is
totally ramified at a prime $\p |p$, $\p \in S -T$. Considering the base change to the strict henselization of $S$ at $\p$ and applying Lemma~\ref{zpext}, we see that each vertical divisor of $\bar X$ in the fibre over $\p$ ramifies in $Kk_\infty$. Replacing $\bar X$ by its normalization in a suitable finite subextension of $Kk_\infty$, we obtain a contradiction using Lemma~\ref{exdiv}.
\end{proof}

Next we consider the case $S=\Spec(\F)$, i.e.\ varieties over a finite field $\F$. In this case we have the degree map
\[
\deg\colon\ \pi_1^t(X)^\ab \lang \pi_1^t(S)^\ab \cong Gal(\bar \F\,|\,\F)\cong \hat \Z,
\]
and we denote the kernel of this degree map by $(\pi_1^t(X)^{\ab})^0$. The image of $\deg$ is an open subgroup of $\hat \Z$ and therefore isomorphic to $\hat \Z$. As $\hat \Z$ is a projective profinite group, we have a (non-canonical) isomorphism
\[
\pi_1^t(X)^\ab\cong (\pi_1^t(X)^{\ab})^0 \times \hat \Z.
\]
Let $p$ be the characteristic of the finite field $\F$. If  $X$ is an open subscheme of a smooth proper variety $\bar X$, then we have a decomposition
\[
(\pi_1^t(X)^{\ab})^0 \cong (\pi_1^\et(X)^{\ab})^0 (\hbox{\rm prime-to-$p$-part}) \oplus (\pi_1^\et(\bar X)^{\ab})^0 (\hbox{\rm $p$-part}),
\]
and both summands are known to be finite. The finiteness statement for $(\pi_1^t(X)^{\ab})^0$ can be generalized to normal schemes.

\begin{theorem}\label{functffinite}
Let $X$ be a normal connected variety over a finite field. Then the group $(\pi_1^t(X)^{\ab})^0$ is finite.
\end{theorem}

\begin{proof} (sketch)
We may replace $X$ by a suitable open subscheme and therefore assume that there exists a smooth morphism $X \lang C$ to a smooth projective curve. Then we proceed in an analogous way as in the proof of Theorem~\ref{finiteness} using the fact that a global field of positive characteristic has exactly one unramified $\hat \Z$-extension, which is obtained by base change from the constant field.
\end{proof}

\section{Tame class field theory}\label{classfieldsec}

In this section  we construct a reciprocity homomorphism from the singular homology group $h_0(X)$ to the abelianized tame fundamental group of an arithmetic scheme~$X$. A sketch of the results of this section is contained in  \cite{S3}.

\ms
Let for the whole section $S=\Spec(\Z)$ and let $X\in \Sch(\Z)$ be connected and regular. If $X$ has $\Real$-valued points, we have to modify the tame fundamental group in the following way:

 We consider the full subcategory of the category of tame coverings of $X$ which consists of that coverings in which every $\Real$-valued point of $X$ splits completely.
After choosing a geometric point $\bar x$ of $X$ we have the fibre functor $(Y\to X)\mapsto \Mor_X(\bar x, Y)$, and its automorphism group $\tilde \pi_1^t(X,\bar x)$ is called the
{\bf modified tame fundamental group} of $X$. It classifies connected tame coverings of $X$ in which every $\Real$-valued point of $X$ splits completely. We have an obvious surjection
\[
 \pi_1^t(X,\bar x) \twoheadrightarrow \tilde \pi_1^t(X,\bar x)
\]
which is an isomorphism if $X(\Real)=\varnothing$.

\ms
For $x\in X(\Real)$ let $\sigma_x\in \pi_1^t(X)^\ab$ be the image of the complex conjugation $\sigma\in \mathit{Gal}(\mathbb C|\Real)$ under the natural map $x_*: \mathit{Gal}(\mathbb C|\Real) \to \pi_1^t(X)^\ab$. By \cite{Sa}, Lemma 4.9 (iii), the map
\[
X(\Real) \lang \pi_1^t(X)^\ab,\ x\longmapsto \sigma_x,
\]
is locally constant for the norm-topology on $X(\Real)$. Therefore the
kernel of the homomorphism
\[
\pi_1^t(X)^\ab \twoheadrightarrow \tilde \pi_1^t(X)^\ab
\]
is an $\F_2$-vector space of dimension less or equal the number of connected components of $X(\Real)$.

\ms
Let $x\in X$ be a closed point. We have a natural isomorphism
\[
\pi_1^\et(\{x\})\cong \mathit{Gal}(\overline{ k(x)}|k(x))\cong\hat \Z,
\]
and we denote the image of the (arithmetic) Frobenius automorphism $\Frob \in G(\overline{ k(x)}|k(x))$
under the natural homomorphism $ \pi_1^\et(\{x\})^{ab} \lang \pi_1^\et(X)^{ab}$
by $\Frob_x$.

\ms
In the following we omit the base scheme $\Spec(\Z)$ from notation, writing
$C_\bullet(X)$ for $C_\bullet(X;\Spec(\Z))$ and similar for homology. Recall that $C_0(X)=Z_0(X)$ is the group of zero-cycles on $X$.  Sending $x$ to $\Frob_x$, we obtain a homomorphism
\[
r: C_0(X) \lang \pi_1(X)^{ab},
\]
which is known to have  dense image (see \cite{La} or \cite{Ras}, lemma
1.7). Our next goal is to show

\begin{theorem} \label{reziex}
The composite map
\[
C_0(X) \mapr{r} \pi_1^\et(X)^{ab} \lang \tilde \pi_1^t(X)^{ab}
\]
factors through $h_0(X)$, thus defining a reciprocity
homomorphism
\[
\rec: h_0(X) \lang \tilde \pi_1^t(X)^{ab},
\]
which has a dense image.
\end{theorem}

In order to prove Theorem~\ref{reziex}, let us apply Theorem~\ref{curve2} to the case of
rings of integers of algebraic number fields. Let $k$ be a finite
extension of $\Q$ and let $\Sigma$ be a finite set of nonarchimedean
primes of $k$. Let $\Cal O_{k,\Sigma}$ be the ring of $\Sigma$-integers of
$k$ and let $E_k^{1,\Sigma}$ be the subgroup of elements in the group of
global units $E_k$ which are $\equiv 1$ at every prime $\p \in \Sigma$.
Let $r_1$ and $r_2$ be the number of real and complex places of $k$. If
$\mathfrak m$ is a product of primes of $k$, then we denote by
$C_{\mathfrak m}(k)$ the ray class group of $k$ with modulus $\mathfrak
m$.

\begin{prop} \label{numbercase1}
For $X= \Spec(\Cal O_{k,\Sigma})$, we have $h_i(X)=0$ for $i\neq 0,1$
and

\ms
 {
 \Item{\rm (i)} $h_0(X)= C_{\mathfrak m}(k)$ with $\mathfrak m= \prod_{\p
 \in \Sigma}\p$.
 \Item{\rm (ii)} $h_1(X)= E_k^{1,\Sigma}\cong \hbox{\rm
 (finite group)} \oplus \Z^{r_1+r_2-1}$.
 }

\ms \noi
In particular, $h_0(X)$ is finite and $h_1(X)$ is
finitely generated. If\/ $\Sigma$ contains at least two primes with
different residue characteristics, then the finite summand in {\rm (ii)}
is zero.
\end{prop}

\begin{proof}
The vanishing of $h_i(X)$ for $i\neq 0,1$ follows
from Theorem~\ref{curve2}. A straightforward computation shows that for
$\mathfrak m= \prod_{\p \in \Sigma}\p$
\[
C_{\mathfrak m}(k)\cong \Pic(\Spec(\Cal O_k), \Sigma),
\]
and the finiteness of $C_{\mathfrak m}(k)$ is well-known. The group
$E_k^{1,\Sigma}$ is of finite index in the full unit group $E_k$.
Therefore the remaining statement in (ii) follows from Dirichlet's unit
theorem. Furthermore, a root of unity, which is congruent to $1$ modulo
two primes of different residue characteristics equals $1$.
\end{proof}

By Theorem~\ref{curve1}, we have an analogous statement for smooth curves over finite fields

\begin{prop} \label{functioncase1} Let $X$ be a smooth, geometrically connected curve over a finite field $\F$ and let $\bar X$ be the uniquely defined smooth compactification of $X$. Let $\Sigma=\bar X-X$ and let $k$ be the function field of $X$. Then we have $h_i(X)=0$ for $i\neq 0,1$
and

\ms
 {
 \Item{\rm (i)} $h_0(X)= C_{\mathfrak m}(k)$ with $\mathfrak m= \prod_{\p  \in \Sigma}\p$.
 \Item{\rm (ii)} $h_1(X)= \left\{ \begin{array}{cl}
                     0& \hbox{ if } \Sigma \neq \varnothing,\\
                     \,\F^\times & \hbox{ if } \Sigma=\varnothing.
                     \end{array} \right.$
}

\ms \noi
In particular, $h_i(X)$ is finite for all $i$.
\end{prop}

\begin{proof}[Proof of Theorem \ref{reziex}] Using Propositions~\ref{numbercase1} and \ref{functioncase1}, classical (one-di\-men\-sio\-nal) class field theory for global fields shows the statement in the case $\dim X=1$. In order to show the general statement,
it suffices by corollary~\ref{erzdD1} to show that for any morphism $f: C \to X$ from a regular curve
$C$ to $X$ and for any $x \in d(C_1(C))$, we have $r(f_*(x))=0$. This
follows from the corresponding result in dimension~$1$ and from the commutative diagram
 \diagram{ccccc}
 d(C_1(C)) & \lang & C_0(C) & \mapr{r_C} & \tilde \pi_1^t(C)^{ab}\\
 \mapd{} & & \mapd{}&& \mapd{}\\
 d(C_1(X)) & \lang & C_0(X) & \mapr{r_X}& \tilde \pi_1^t(X)^{ab}.
 \enddiagram
\end{proof}

\medskip\noindent
In order to investigate the reciprocity map, we use Wiesend's version of higher dimensional class field theory \cite{W2}. We start with the arithmetic case, i.e.\ when $X$ is flat over $\Spec(\Z)$. In this case  $\tilde \pi_1^t(X)^\ab$ is finite by Theorem~\ref{finiteness}.

\begin{theorem}\label{tameclassfield}
Let $X$ be a regular, connected scheme, flat and of finite type over $\Spec(\Z)$. Then the reciprocity homomorphism
\[
\rec_X\colon h_0(X) \lang \tilde \pi_1^t(X)^\ab
\]
is an isomorphism of finite abelian groups.
\end{theorem}

\begin{remark}
If $X$ is proper, then $h_0(X)\cong \CH_0(X)$ and $\tilde \pi_1^t(X)^\ab\cong \tilde \pi_1^\et(X)^\ab$, and we recover the unramified class field theory for arithmetic schemes of Bloch and Kato/Saito \cite{K-S}, \cite{Sa}.
\end{remark}

\begin{proof}[Proof of Theorem~\ref{tameclassfield}] Recall the definition of Wiesend's id\`{e}le group ${\cal I}_X$. It is defined by
\[
{\cal I}_X:=  Z_0(X) \oplus \bigoplus_{C\subset X}\bigoplus_{v\in C_\infty}k(C)_v^\times.
\]
Here $C$ runs through all closed integral subschemes of $X$ of dimension~$1$, $C_\infty$ is the finite set of places (including the archimedean ones if $C$ is horizontal) of the global field $k(C)$ with center outside $C$ and $k(C)_v$ is the completion of $k(C)$ with respect to $v$.
${\cal I}_X$ becomes a topological group by endowing the group $Z_0(X)$ of zero cycles on $X$ with the discrete topology, the groups $k(C)_v^\times$ with their natural locally compact topology and the direct sum with the direct sum topology.\footnote{The topology of a finite direct sum is just the product topology, and the topology of an infinite direct sum is the direct limit topology of the finite partial sums.}

\medskip
The id\`{e}le class group ${\cal C}_X$ is defined as the cokernel of the natural map
\[
\bigoplus_{C\subset X}k(C)^\times \lang  {\cal I}_X.
\]
which is given for a fixed $C\subset X$ by the divisor map $k(C)^\times \to Z_0(C) \to Z_0(X)$ and the diagonal map $k(C)^\times \to \bigoplus_{v\in C_\infty}k(C)_v^\times$. ${\cal C}_X$ is endowed the quotient topology of ${\cal I}_X$.

We consider the quotient ${\cal C}_X^t$ of ${\cal C}_X$ obtained by cutting out the $1$-unit groups at all places outside $X$. More precisely,
let for $v\in C_\infty$, $U^1(k(C)_v)$ be the group of principal units in the local field $k(C)_v$. We make the notational convention  $U^1(K)=K^\times$ for the archimedean local fields $K=\Real,\mathbb C$. Then
\[
{\cal U}^t_X:= \bigoplus_{C\subset X}\bigoplus_{v\in C_\infty}U^1(k(C)_v)
\]
is an open subgroup of the id\`{e}le group ${\cal I}_X$ and we put
\[
{\cal C}_X^t:= \coker(\bigoplus_{C\subset X}k(C)^\times \lang  {\cal I}_X/{\cal U}^t_X ).
\]
Consider the map
\[
R:\ {\cal I}_X \lang \pi_1^\et(X)^\ab
\]
which is given by the map $r: Z_0(X)\to \pi_1(X)^\ab$ defined above and the reciprocity maps of local class field theory
\[
\rho_v:\  k(C)_v^\times \lang \pi_1^\et(\Spec(k(C)_v))^\ab
\]
followed by the natural maps $\pi_1^\et(\Spec(k(C)_v))^\ab \to \pi_1^\et(X)^\ab$ for all $C\subset X$, $v\in C_\infty$.
By the \cite{W2} Theorem~1\,(a) the homomorphism $R$ induces an isomorphism
\[
\rho:\ {\cal C}_X^t \liso \tilde \pi_1^t(X)^\ab.
\]
Now we consider the obvious map
\[
\phi:\ Z_0(X) \lang {\cal C}_X^t.
\]
The kernel of $\phi$ is the subgroup in $Z_0(X)$ generated by elements of the form $\div(f)$ where $C\subset X$ is a closed curve and  $f$ is an invertible rational function on $C$ which is in $U^1(k(C)_v)$ for all $v\in C_\infty$. By Theorem~\ref{h0char} we obtain $\ker(\phi)=d_1(C_1(X))$. Therefore $\phi$ induces an injective homomorphism
\[
i:\ h_0(X) \hookrightarrow {\cal C}_X^t
\]
with $\rho \circ i = \rec$. As $\rho$ is injective, $\rec$ is injective, and hence an isomorphism.
\end{proof}
Finally, assume that $\bar X$ is  regular, flat and proper over $\Spec(\Z)$, let $D\subset X$ be a divisor and $X=\bar X -D$. In \cite{S2} we introduced the relative Chow group of zero cycles $\CH_0(\bar X,D)$ and constructed, under a mild technical assumption, a reciprocity isomorphism $\rec': \CH_0(\bar X,D)\stackrel{\sim}{\to}\tilde \pi_1^t(X)^\ab$. By \cite{S2}, Proposition~2.4, there exists natural projection $\pi: h_0(X)\twoheadrightarrow \CH_0(\bar X,D)$ with $\rec=\rec'\circ \pi$. We obtain the

\begin{theorem}
Let $\bar X$ be a regular, connected scheme, flat and proper over $\Spec(\Z)$, such that its
generic fibre $X\otimes_\Z \Q$ is projective over $\Q$. Let $D$ be a divisor on
$\bar X$ whose vertical irreducible components are normal schemes. Put $X=\bar X-D$. Then the natural homomorphism
\[
h_0(X) \lang \CH_0(\bar X,D)
\]
is an isomorphism of finite abelian groups.
\end{theorem}

Finally, we deal with the geometric case. The next theorem was proved in 1999 by M.~Spie{\ss} and the author under the assumption that $X$ has a smooth projective compactification, see \cite{S-S1}.  Now we get rid of this assumption.

\begin{theorem}\label{s-spiess}
Let $X$ be a smooth, connected variety over a finite field $\F$. Then the reciprocity homomorphism
\[
\rec_X\colon h_0(X) \lang \pi_1^t(X)^\ab
\]
is injective. The image of\/ $\rec_X$ consists of all elements whose degree in $\mathit{Gal}(\bar \F |\F)$ is an integral power of the Frobenius automorphism.  In particular, the cokernel $\coker(\rec_X)\cong \hat\Z/\Z$ is uniquely divisible. The induced map on the degree-zero parts $\rec_X^0\colon h_0(X)^0 \stackrel{\sim}{\to} (\pi_1^t(X)^\ab)^0$ is an isomorphism of finite abelian groups.
\end{theorem}

\begin{proof} The proof is strictly parallel to the proof of Theorem~\ref{tameclassfield}, using Theorem~\ref{h0char} and the tame version of Wiesend's class field theory for smooth varieties over finite fields \cite{W2}, Theorem~1\,(b).
\end{proof}

\vskip2cm

\noindent \footnotesize{Alexander Schmidt, NWF I - Mathematik, Universit\"{a}t Regensburg, D-93040
Regensburg, Deutschland. email: alexander.schmidt@mathematik.uni-regensburg.de}
\end{document}